\documentclass[twoside,reqno]{amsart}

\usepackage{framed}

\usepackage{amsmath}
\usepackage{amssymb}
\usepackage{amsthm}
\usepackage{mathtools}
\usepackage{graphicx}
\usepackage{amsfonts}

\usepackage{accents}

\makeatletter
\def\widebar{\accentset{{\cc@style\underline{\mskip10mu}}}}
\makeatother

\usepackage{array}

\newcolumntype{L}[1]{>{\raggedright\arraybackslash}m{#1}}
\newcolumntype{C}[1]{>{\centering\arraybackslash}m{#1}}

\usepackage[table]{xcolor}

\usepackage[colorlinks=true,linkcolor=red]{hyperref}
\newtheorem{theorem}{Theorem}[section]
\newtheorem{proposition}[theorem]{Proposition}
\newtheorem{lemma}[theorem]{Lemma}
\newtheorem{corollary}[theorem]{Corollary}

\theoremstyle{definition}
\newtheorem{definition}[theorem]{Definition}
\newtheorem{example}[theorem]{Example}
\newtheorem{remark}[theorem]{Remark}

\numberwithin{equation}{section}

\DeclareSymbolFont{stmry}{U}{stmry}{m}{n}
\DeclareMathDelimiter\llbracket{\mathopen}{stmry}{"4A}{stmry}{"71}
\DeclareMathDelimiter\rrbracket{\mathclose}{stmry}{"4B}{stmry}{"79}

\newcommand{\CH}{\mathrm{C}_{\mkern-1.2mu\mathfrak H}\mkern-1.9mu}
\newcommand{\llf}{\llbracket}
\newcommand{\rlf}{\rrbracket}
\newcommand{\bllf}{\bigl\llbracket}
\newcommand{\brlf}{\bigr\rrbracket}

\newcommand{\lda}{\langle\mkern-4.5mu\langle}
\newcommand{\rda}{\rangle\mkern-4.5mu\rangle}

\newcommand{\kip}{[{\mkern1.7mu{\cdot}\mkern0.9mu,\mkern-2.1mu{\cdot}\mkern1.4mu}]}

\newcommand{\ahip}{%
\langle{\mkern1.6mu{\cdot}\mkern0.9mu,\mkern-2.1mu{\cdot}\mkern1.6mu}\rangle}
\newcommand{\nor}{{\|\mkern1mu{\cdot}\mkern1mu\|}}
\newcommand{\lfkip}{%
\llf{\mkern1.6mu{\cdot}\mkern0.9mu,\mkern-2.1mu{\cdot}\mkern1.6mu}\rlf}
\newcommand{\dahip}{%
\lda{\mkern1.6mu{\cdot}\mkern0.9mu,\mkern-2.1mu{\cdot}\mkern1.6mu}\rda}

\newcommand{\iu}{\mkern0.35mu\mathrm{i}\mkern0.55mu}

\def\-{\raisebox{.05pt}{-}}

\newcommand{\gt}{>}

 \DeclareMathOperator{\lspan}{span}

 \DeclareMathOperator{\ran}{ran}
 \DeclareMathOperator{\dom}{dom}
 \DeclareMathOperator{\nul}{nul}

\begingroup
\def\2#1{\ifnum#1<10 0\fi\the#1}
\xdef\isodayandtime%
{\the\year-\2\month-\2\day\space\2{\count0}:\2{\count2}}

\endgroup

\allowdisplaybreaks

\makeatletter \renewcommand{\p@enumii}{} \makeatother
\makeatletter \renewcommand{\p@enumiii}{} \makeatother

\begin{document}

\title[Compressions in Krein Spaces]{ %
Compressions of Linear Relations in Krein Spaces}

\author{B.~\'{C}urgus}
\address{Department of Mathematics, Western Washington University, Bellingham, WA 98225, USA}
\email{\texttt{curgus@wwu.edu}}
\author{A.~Dijksma}
\address{Bernoulli Institute of Mathematics, Computer Science and Artificial Intelligence \\ University of Groningen \\
P.O. Box 407 \\
9700 AK Groningen, The Netherlands}
\email{\texttt{a.dijksma@rug.nl}}

\begin{abstract}
Given an orthogonal decomposition \(\widetilde{\mathfrak{H}} = \mathfrak{H}[\dotplus]\widehat{\mathfrak{H}}\) of a Krein space, we study compressions \(\CH(R)\) of linear relations \(R\) in \(\widetilde{\mathfrak{H}}\) to \(\mathfrak{H}\). We formulate a sufficient condition that ensures that compression commutes with taking adjoints, that is, the adjoint identity \({\mathrm{C}_{\mkern-1.2mu\mathfrak H}\mkern-1.9mu}(R^{\widetilde{*}}) = {\mathrm{C}_{\mkern-1.2mu\mathfrak H}\mkern-1.9mu}(R)^{*}\) holds. We derive preservation and reflection results for the compression of maximal dissipative, maximal symmetric, and self-adjoint linear  relations. For nonnegative linear relations, the adjoint identity also implies that compression commutes with taking the Friedrichs extension. For \(\mathfrak{H}\ne\{0\}\), we characterize finite dimensionality of the exit space \(\widehat{\mathfrak{H}}\) by universal compression properties. Finally, we give a criterion for the adjoint identity to hold in terms of bounded linear functionals on \(R\).
\end{abstract}



\maketitle

\tableofcontents

\section{Introduction}
In this paper we study compressions of linear relations in Krein spaces. All vector spaces in this paper are complex. We use standard notation from Krein spaces, as in \cite{AzB}, see also \cite{Bog}, and the theory of linear relations therein, see \cite{DS}. The basic definitions and statements are recalled in Section~\ref{sec-Prel}. We consider a  Krein space \(\bigl(\widetilde{\mathfrak H},\kip_{\widetilde{\mathfrak H}} \bigr)\) and assume that it is the orthogonal direct sum of its Krein subspaces \(\mathfrak H\) and \(\widehat{\mathfrak H}\): \(\widetilde{\mathfrak{H}} = \mathfrak{H} [\dot{+}] \widehat{\mathfrak H}\). The orthogonal projections in \(\widetilde{\mathfrak H}\) onto these subspaces are denoted by \(\widetilde{P}_{\mathfrak H}\) and \(\widetilde{P}_{\widehat{\mathfrak H}}\). We denote adjoints of linear relations in Krein spaces \(\bigl(\widetilde{\mathfrak{H}},\kip_{\widetilde{\mathfrak{H}}}\bigr)\)
and \(\bigl(\mathfrak{H},\kip_{\mathfrak{H}}\bigr)\) by \(^{\widetilde{*}}\) and \(^{*}\), respectively. Here \(\kip_{\mathfrak{H}}\) is the restriction of \(\kip_{\widetilde{\mathfrak{H}}}\) to \(\mathfrak{H}^2\). 

Let \(R\) be a linear relation in \(\widetilde{\mathfrak H}\). The \emph{compression of} \(R\) \emph{to} \(\mathfrak{H}\) is the linear relation in \(\mathfrak{H}\) defined by 
\[
\CH(R) = \Bigl\{ \{f, \widetilde{P}_{\mathfrak H} \widetilde{g}\} : 
f \in \mathfrak{H}, \ \{f, \widetilde{g}\} \in R \Bigr\} =\widetilde{P}_{\mathfrak H}\bigl( R \cap \bigl( \mathfrak{H} \times  \widetilde{\mathfrak{H}}\bigr) \bigr). 
\]
In \cite[p.~25]{WeSt} it is called the part of \(R\) in \(\mathfrak H\). It is the main object of interest in this paper. In this definition \(\mathfrak H\) and \(\widehat{\mathfrak H}\) play different roles, \(\mathfrak H\) is called the \emph{base space} and \(\widehat{\mathfrak H}\) the \emph{exit space}.

In 1968 Stenger \cite{St1968} proved the following result.
\begin{theorem} \label{Stenger}
Let \(\widetilde{\mathfrak{H}} = \mathfrak{H} [\dot{+}] \widehat{\mathfrak H}\) be a Hilbert space and assume that \(\widehat{\mathfrak H}\) is finite dimensional. If \(R\) is a densely defined self-adjoint operator in \(\widetilde{\mathfrak H}\), then \(\CH(R)\) is a densely defined self-adjoint operator in \(\mathfrak H\).
\end{theorem}

In 2011 this theorem was generalized by Nudelman \cite{Nu} to maximal dissipative operators.
\begin{theorem} \label{th-nud}
Let \(\widetilde{\mathfrak{H}} = \mathfrak{H} [\dot{+}] \widehat{\mathfrak H}\) be a Hilbert space and assume that \(\widehat{\mathfrak H}\) is finite dimensional. If \(R\) is a densely defined maximal dissipative operator in \(\widetilde{\mathfrak H}\), then \(\CH(R)\) is a densely defined maximal dissipative operator in \(\mathfrak H\).
\end{theorem}
In 2012 and 2013 these two statements were generalized to Krein spaces, linear relations and ``if and only if'' statements, see \cite{AD2012} and \cite{ADW2013}.

\begin{theorem} \label{th-mdiff}
Let \(\widetilde{\mathfrak H} = \mathfrak{H} [\dot{+}] \widehat{\mathfrak H}\) be a Krein space and assume that \(\widehat{\mathfrak H}\) is finite dimensional. Let \(R\) be a closed dissipative linear relation in \(\widetilde{\mathfrak H}\). Then \(R\) is maximal dissipative if and only if \(R(0)=R^{\widetilde{*}}\mkern-2mu(0)\) and \(\CH(R)\) is a maximal dissipative linear relation in \(\mathfrak H\).
\end{theorem}

\begin{theorem} \label{th-msiff}
Let \(\widetilde{\mathfrak H} = \mathfrak{H} [\dot{+}] \widehat{\mathfrak H}\) be a Krein space and assume that \(\widehat{\mathfrak H}\) is finite dimensional. Let \(R\) be a closed symmetric linear relation in \(\widetilde{\mathfrak H}\). Then \(R\) is maximal symmetric (respectively, self-adjoint) if and only if \(R(0)=R^{\widetilde{*}}\mkern-2mu(0)\) and \(\CH(R)\) is a maximal symmetric (respectively, self-adjoint) linear relation in \(\mathfrak H\).
\end{theorem}

The maximal dissipative and self-adjoint statements are
\cite[Theorems~3.3 and~4.1]{ADW2013}, respectively.
The maximal symmetric statement follows from the former, see Theorem~\ref{th-maxnnR}.

In this paper we show that the finite dimensionality of \(\widehat{\mathfrak{H}}\) is exactly the condition under which several properties of linear relations \(R\) in \(\widetilde{\mathfrak{H}}\) are universal, that is, properties that hold for all \(R\) in \(\widetilde{\mathfrak{H}}\) or all closed \(R\) in \(\widetilde{\mathfrak{H}}\). We prove, in Section~\ref{SecEqi}, the following theorem.

\begin{theorem}\label{equivalences}
Let \(\widetilde{\mathfrak H} = \mathfrak{H} [\dot{+}] \widehat{\mathfrak H}\) be a Krein space with \(\mathfrak{H}\neq \{0\}\). The following statements are equivalent.
\begin{enumerate}
\renewcommand*\theenumi{\roman{enumi}}
\renewcommand*\labelenumi{{\rm (\theenumi)}}
\item \label{equivalences-i1}
\(\widehat{\mathfrak H}\) is finite dimensional.
\item  \label{equivalences-i1a}
For every closed linear relation \(R\) in \(\widetilde{\mathfrak H}\) the subspace \(\widetilde{P}_{\widehat{\mathfrak H}}(\dom R)\) is closed. 
\item  \label{equivalences-i1b}
For every closed linear relation \(R\) in \(\widetilde{\mathfrak H}\) the linear relation \(\widetilde{P}_{\mathfrak H} R\)  in \(\widetilde{\mathfrak H}\) is closed.
\item \label{equivalences-cc}
For every linear relation \(R\) in \(\widetilde{\mathfrak H}\) we have \(\CH(\widebar R)=\overline{\CH(R)}\).
\item  \label{equivalences-i20}
For every linear relation \(R\) in \(\widetilde{\mathfrak H}\) we have \(\CH(R^{\widetilde{*}})=\CH(R)^*\).
\item  \label{equivalences-i2}
For every closed linear relation \(R\) in \(\widetilde{\mathfrak H}\) we have \(\CH(R^{\widetilde{*}})=\CH(R)^*\).
\item \label{equivalences-i3}
For every maximal dissipative linear relation \(R\) in \(\widetilde{\mathfrak H}\)  its compression  \(\CH(R)\) is a maximal dissipative linear relation in \(\mathfrak H\).
\item \label{equivalences-i4}
For every maximal symmetric linear relation \(R\) in \(\widetilde{\mathfrak H}\)  its compression \(\CH(R)\) is a maximal symmetric linear relation in \(\mathfrak H\).
\item \label{equivalences-i5}
For every self-adjoint linear relation \(R\) in \(\widetilde{\mathfrak H}\) its compression  \(\CH(R)\) is a self-adjoint linear relation in \(\mathfrak{H}\).
\item \label{equivalences-Friedrichs} 
For every nonnegative linear relation \(S\) in \(\widetilde{\mathfrak H}\) we have that \(\CH(S)\) is a nonnegative linear relation in \(\mathfrak{H}\) and \(\CH(S_{\mathrm F}) = (\CH(S))_{\mathrm F}\), where the subscript \({\mathrm F}\) stands for the Friedrichs extension of a nonnegative linear relation. 
\end{enumerate}
\end{theorem}

These statements are reminiscent of the theorems by Gustafson \cite{Gus} in 1969 and Van Casteren and Goldberg \cite{vCGo} in 1970 which we present in one statement.

\begin{theorem}
Let \(S\) be a closed densely defined linear operator in a Hilbert space \(\mathfrak K\). Then \((TS)^* = S^*T^*\) holds for all closed densely defined linear operators \(T\) on \(\mathfrak K\) if and only if \(\ran S\) has finite codimension.
\end{theorem}

The ``if'' statement is due to Gustafson (for a Banach space version see Schechter \cite{Sch1970}) and can be used to prove Theorem~\ref{Stenger}. The ``only if''  statement is due to van Casteren and Goldberg. The proof consists of constructing a \(T\) for which the equality of adjoints does not hold if the codimension of the range \(\ran S\) is infinite. Similarly, in the proof of Theorem~\ref{equivalences} we also construct counterexamples.

We conclude this introduction by briefly outlining the contents of each section. Section~\ref{sec-Prel} establishes preliminaries regarding subspaces and linear relations in Krein spaces, while Section~\ref{sec-comps} presents basic tools related to compressions in Krein spaces.

In Section~\ref{SecCAACn} we consider closed linear relations \(R\) and give sufficient conditions for compression to commute with taking adjoints. One of these conditions, namely the closedness of the linear relation  \(\widetilde{P}_{\mathfrak{H}} R^{\widetilde{*}}\), occurs as a condition in \cite[(4.3) and Lemma~4.4]{ADW2013}. Theorem~\ref{th-domRandPR*} shows that this condition is equivalent to the closedness of the subspace 
\(\widetilde{P}_{\widehat{\mathfrak{H}}} (\operatorname{dom}R)\) in \(\widehat{\mathfrak{H}}\). Theorem~\ref{th-comadjn} shows that these equivalent
conditions imply the adjoint identity \(\CH(R^{\widetilde{*}})=\CH(R)^*\). It follows that they hold automatically when the exit space \(\widehat{\mathfrak{H}}\) is finite dimensional, see also Theorem~\ref{equivalences}.

We organize Section~\ref{SecPR} in two complementary directions. Firstly, under the assumption that
\(\widetilde P_{\widehat{\mathfrak H}}\bigl(\dom R^{\widetilde{*}}\bigr)\) is closed in \(\widehat{\mathfrak H}\), Theorem~\ref{th-comadj4} gives preservation results for compressions. 
By this we mean that certain properties of \(R\) are inherited by its compression   \(\CH(R)\).   
Theorem~\ref{th-comadj4} generalizes Theorems~\ref{Stenger} and~\ref{th-nud}. For other preservation results see \cite{ACD16}. Secondly, under the assumption that \(\widetilde P_{\widehat{\mathfrak H}}(\dom R)\) is closed in \(\widehat{\mathfrak H}\), Theorem~\ref{th-conv} gives reflection results. By this we mean that with additional conditions, the properties of the compression \(\CH(R)\) imply the same properties for \(R\). Combining these two directions yields the four ``if and only if'' statements in Theorem~\ref{th-4-iff}, which generalize Theorems~\ref{th-mdiff} and~\ref{th-msiff}.

In Section~\ref{sec-Friedrichs} we prove that  \(\CH(S^{\widetilde{*}})=\CH(S)^*\) implies \(\CH(S_{\mathrm F})=(\CH(S))_{\mathrm F}\) for every nonnegative linear relation \(S\) in \(\widetilde{\mathfrak{H}}\). This extends \cite[Theorem~4.5]{ADW2013} from closed linear relations in Hilbert spaces to arbitrary nonnegative  linear relations in Krein spaces and permits infinite-dimensional exit spaces. When the exit space is finite dimensional, the adjoint identity holds automatically. Thus the hypothesis that at least one of the conditions in \cite[(4.3)]{ADW2013} holds is unnecessary for the Friedrichs compression identity \(\CH(S_{\mathrm F})=(\CH(S))_{\mathrm F}\). The definition of the Friedrichs extension of a nonnegative linear relation in a Krein space is analogous to the definition in a Hilbert space. For background on self-adjoint extensions of nonnegative operators in Krein spaces, see \cite{Curgus1989}.

In Section~\ref{SecEqi} we give a proof of Theorem~\ref{equivalences}. 

In Appendix~\ref{Sec-funs} we present a necessary and sufficient boundedness criterion for compression to commute with taking adjoints and we compare it with the closedness conditions used in Section~\ref{SecCAACn}.

\section{Krein space preliminaries} \label{sec-Prel}

\subsection{Subspaces of a Krein space} \label{subs-Ks}
Our terminology follows \cite{AzB} with the exception of the concept of a subspace. By a \emph{subspace} of a topological vector space we mean a subspace of a vector space  disregarding the topology. If a subspace is a closed subset of a topological space we will say so, and call it a \emph{closed subspace}.

Let \(\bigl(\mathfrak K,\kip_{\mathfrak K}\bigr)\) be a Krein space and let
\(\mathfrak A\) be a subset of \(\mathfrak K\). We say that \(\mathfrak A\) is
\emph{nonnegative}, \emph{nonpositive}, or \emph{neutral} if, respectively,
\[
[x,x]_{\mathfrak K}\geq 0,\qquad
[x,x]_{\mathfrak K}\leq 0,\quad \text{or} \quad
[x,x]_{\mathfrak K}=0
\]
for all \(x\in\mathfrak A\). The subset \(\mathfrak A\) is called
\emph{semidefinite} if it is either nonnegative or nonpositive. The
\emph{orthogonal complement} of \(\mathfrak A\), denoted by
\(\mathfrak A^{[\perp]}\), is defined by
\[
\mathfrak A^{[\perp]}
=
\bigl\{x\in\mathfrak K:[x,a]_{\mathfrak K} = 0  \ \
\text{for all} \ \ a\in\mathfrak A \bigr\}.
\]
For a subspace \(\mathfrak{L}\) of \(\mathfrak K\) the subspace \(\mathfrak{L} \cap \mathfrak{L}^{[\perp]}\) is the \emph{isotropic part} of \(\mathfrak{L}\). 

A nonnegative subspace \(\mathfrak{L}\) of \(\mathfrak K\) is called \emph{maximal nonnegative} if, whenever \(\mathfrak M\) is a nonnegative subspace of \(\mathfrak K\) such that
\(\mathfrak L\subseteq\mathfrak M\),
then \(\mathfrak L=\mathfrak M\). Maximal nonpositive and maximal neutral subspaces are defined analogously. A subspace \(\mathfrak L\) of \(\mathfrak K\) is called \emph{hyper-maximal neutral} if it is both maximal nonnegative and maximal nonpositive. 

Since the closure of a nonnegative, nonpositive, or neutral subspace retains the same property, it follows that every maximal nonnegative, nonpositive, or neutral subspace is closed.

The following characterizations are in essence proved in \cite{AzB}. In the proof we give exact references.

\begin{theorem} \label{th-maxsub}
Let \((\mathfrak{K},\kip_\mathfrak{K})\) be a Krein space and let \(\mathfrak{L}\) be a subspace of \(\mathfrak{K}\). The following statements hold.
\begin{enumerate}
\renewcommand*\theenumi{\alph{enumi}}
\renewcommand*\labelenumi{{\rm (\theenumi)}}
\item  \label{th-pr-maxsubs-i1}
\(\mathfrak L\) is maximal nonnegative (respectively, maximal nonpositive)
in \(\mathfrak K\) if and only if \(\mathfrak L\) is closed and nonnegative
(respectively, nonpositive) and \(\mathfrak L^{[\perp]}\) is nonpositive
(respectively, nonnegative) in \(\mathfrak K\).

\item \label{th-pr-maxsubs-i2}
\(\mathfrak L\) is maximal neutral in \(\mathfrak K\) if and only if \(\mathfrak L\) is neutral and is maximal nonnegative or maximal nonpositive in \(\mathfrak K\).

\item  \label{th-pr-maxsubs-i3}
\(\mathfrak{L}\) is hyper-maximal neutral in \(\mathfrak{K}\) if and only if \(\mathfrak{L}=\mathfrak{L}^{[\perp]}\).
\end{enumerate}
\end{theorem}

\begin{proof}
(\ref{th-pr-maxsubs-i1})
The ``only if'' part is  \cite[Theorem~1.8.11]{AzB}. To prove the ``if'' part assume that \(\mathfrak{L}\) is closed and nonnegative and \(\mathfrak{L}^{[\perp]}\) is nonpositive. By \cite[Remark~1.8.8]{AzB} there exists a maximal nonnegative subspace \(\mathfrak{M}\) such that \(\mathfrak{L} \subseteq \mathfrak{M}\). Then, by \cite[Theorem~1.8.11]{AzB} \(\mathfrak{M}^{[\perp]}\) is maximal nonpositive. Since  \(\mathfrak{M}^{[\perp]}\subseteq \mathfrak{L}^{[\perp]}\) and  \(\mathfrak{L}^{[\perp]}\) is also nonpositive, this yields \(\mathfrak{L}^{[\perp]} = \mathfrak{M}^{[\perp]}\). Using that \(\mathfrak{L}\) and \(\mathfrak{M}\) are closed and taking the orthogonal complement again we obtain \(\mathfrak{L}= \mathfrak{M}\). This proves that \(\mathfrak{L}\) is maximal nonnegative. The proof for the maximal nonpositive case is similar.

(\ref{th-pr-maxsubs-i2})
The ``only if'' part is \cite[Corollary~1.4.14]{AzB}. The ``if'' part follows from (\ref{th-pr-maxsubs-i1})  and \cite[Proposition~1.4.18]{AzB}.

(\ref{th-pr-maxsubs-i3})
This is \cite[Proposition~1.4.19]{AzB}.
\end{proof}

\subsection{Linear relations in a Krein space} \label{subs-ddip}
Let \(\bigl(\mathfrak{K},\kip_\mathfrak{K} \bigr)\) be a Krein space and let \(J\) be a fundamental symmetry on it; that is, \(J: \mathfrak K \to \mathfrak K\) is a bounded linear operator, \(J^2 = I_{\mathfrak K}\) and \(\mathfrak{K}\) equipped with the inner product \(\ahip_\mathfrak{K} = [J\mkern+1.2mu \cdot\mkern+1mu,\mkern-2.5mu\cdot\mkern+2mu ]_\mathfrak{K}\) is a Hilbert space.

A \emph{linear relation} in \(\mathfrak{K}\) is a subspace of the Cartesian product \(\mathfrak{K}^2\). Every linear operator in \(\mathfrak K\) may be viewed as a linear relation by identifying it with its graph. For the definitions of the domain \(\dom R\), the range \(\ran R\), the null space \(\nul R\), the inverse \(R^{-1}\), and the usual linear operations with \(R\), we refer to \cite{DS} and \cite[Chapter~2]{AzB}. We often use the multivalued part of \(R\) which we denote by \(R(0)\):
\[
R(0) = \bigl\{y \in \mathfrak{K} :  \{0,y\} \in R \bigr\}. 
\]

To recall the adjoint of \(R\) and describe it geometrically, it is convenient to equip \(\mathfrak{K}^2\) with the \emph{Lagrange indefinite inner product}:
\begin{equation} \label{bdHf}
\bllf \{f,g\}, \{h,k\} \brlf_{\mathfrak{K}^2} =
\iu\bigl( [ f, k ]_\mathfrak{K} - [g,h ]_\mathfrak{K} \bigr), \qquad \{f,g\}, \{h,k\} \in \mathfrak{K}^2.
\end{equation}
Define \(\mathbb{J}: \mathfrak{K}^2 \to \mathfrak{K}^2\) by
\begin{equation} \label{eq-dbbJ}
\mathbb{J}\{f,g\} =  \bigl\{- \iu J g, \iu J f \bigr\},
\qquad \{f,g\}  \in \mathfrak{K}^2.
\end{equation}
Then \(\mathbb J^2 = I_{\mathfrak K^2}\), and for all \( \{f,g\}, \{h,k\}  \in \mathfrak{K}^2\) we have
\begin{align*}
\bllf \mathbb J\{f,g\},\{h,k\}\brlf_{\mathfrak{K}^2}
&=
\iu\bigl([-\iu Jg,k]_{\mathfrak K}-[\iu Jf,h]_{\mathfrak K}\bigr)\\
&=
\langle f,h\rangle_{\mathfrak K}
+
\langle g,k\rangle_{\mathfrak K}.
\end{align*}
Thus, \(\mathfrak{K}^2\) with the inner product \(\llf \mathbb{J} \mkern2mu \cdot\mkern+1mu,\mkern-2.5mu\cdot\mkern+2mu \rlf_{\mathfrak{K}^2}\) is the Hilbert space direct product of two copies of the Hilbert space \(\bigl(\mathfrak{K}, \ahip_\mathfrak{K}\bigr)\). Therefore we write \(\llf \mathbb{J} \mkern2mu \cdot\mkern+1mu,\mkern-2.5mu\cdot\mkern+2mu \rlf_{\mathfrak{K}^2} = \dahip_{\mathfrak{K}^2}\). Since \(\bigl(\mathfrak{K}^2,\dahip_{\mathfrak{K}^2}\bigr)\) is a Hilbert space, \(\mathbb J\) is a fundamental symmetry on \(\bigl(\mathfrak K^2,\lfkip_{\mathfrak{K}^2}\bigr)\). This proves that \(\bigl(\mathfrak K^2,\lfkip_{\mathfrak{K}^2}\bigr)\) is a Krein space. It is called the \emph{space of graphs} or the \emph{Lagrange-Krein space} associated with the Krein space \(\bigl(\mathfrak{K},\kip_\mathfrak{K} \bigr)\).

The adjoint of a linear relation \(R\) in \(\bigl(\mathfrak{K},\kip_\mathfrak{K} \bigr)\), denoted by \(R^{[*]}\),  is precisely the orthogonal complement of the subspace \(R\) in the Krein space \(\bigl(\mathfrak{K}^2,\lfkip_{\mathfrak{K}^2}\bigr)\):
\[
R^{[*]} = R^{\llf\perp\rlf}= \Bigl\{
\{x,y\} \in \mathfrak{K}^2 : [f, y]_{\mathfrak{K}} - [g,x]_{\mathfrak{K}} = 0 \ \ \text{for all} \ \ \{f,g\} \in R
\Bigr\}.
\]

For later use, note that for a linear relation \(R\) in \(\bigl(\mathfrak{K},\kip_\mathfrak{K} \bigr)\) we have 
\begin{equation} \label{eq-R*(0)} 
R^{[*]}\mkern-2.5mu(0) = (\dom R)^{[\perp]}.
\end{equation} 
This equality follows from the equivalences
\begin{equation*}
x \in (\dom R)^{[\perp]}  \quad \Leftrightarrow \quad 
\forall\mkern2mu \{f,g\} \in R \ \ \ \llf \{f,g\}, \{0,x\} \rlf = 0  \quad \Leftrightarrow \quad
\{0,x\} \in R^{[*]}. 
\end{equation*}

The standard concepts related to linear relations have the geometric analogs presented in the table below. The table provides a dictionary and, at the same time,  the definitions of the listed concepts for linear relations.
\begin{center}
\renewcommand{\arraystretch}{1.4}

\setlength{\tabcolsep}{10pt}

\begin{tabular}{|L{0.51\textwidth}|L{0.37\textwidth}|}
\hline
\rowcolor{gray!10}
\multicolumn{1}{|C{0.51\textwidth}|}{%
\shortstack[c]{\rule{0pt}{10pt}Property of a linear relation \(R\)\\
in \(\bigl(\mathfrak K,\kip_{\mathfrak K}\bigr)\)}}
&
\multicolumn{1}{C{0.37\textwidth}|}{%
\shortstack[c]{\rule{0pt}{10pt}Property of \(R\) as a subspace\\
of \(\bigl(\mathfrak K^2,\lfkip_{\mathfrak{K}^2}\bigr)\)}}
\\[8pt]
\hline%\hline
adjoint relation: \ \ \(R^{[*]}\)
& orthogonal complement
\\[2pt]\hline
dissipative: \ \
\(\operatorname{Im}\mkern0.5mu [g,f]_{\mathfrak K} \geq 0
\mkern15mu \forall \mkern1.5mu \{f,g\} \in R\)
& nonnegative
\\[2pt]\hline
symmetric: \ \ \(R \subseteq R^{[*]}\)
& neutral
\\[2pt]\hline
maximal dissipative
& maximal nonnegative
\\[2pt]\hline
maximal symmetric
& maximal neutral
\\[2pt]\hline
self-adjoint: \ \ \(R = R^{[*]}\)
& hyper-maximal neutral
\\[2pt]\hline
\end{tabular}
\end{center}
For example, the second row of the above table defines a linear relation \(R \subseteq \mathfrak{K}^2\) to be \emph{dissipative} in \(\bigl(\mathfrak K,\kip_{\mathfrak K}\bigr)\) if and only if for all  \(\{f,g\} \in R\) we have
\begin{equation} \label{eq-dissym}
\operatorname{Im}\mkern0.5mu [g,f]_{\mathfrak K}=
\frac{1}{2}\iu \bigl([f,g]_{\mathfrak K} - [g,f]_{\mathfrak K}\bigr)
=\frac{1}{2}
\bigl\llf \{f,g\}, \{f,g\} \bigr\rlf_{\mathfrak{K}^2} \geq 0.
\end{equation}

Since all vector spaces in this paper are complex, the polarization identity applied to the inner product space \(\bigl(R, \lfkip_{\mathfrak{K}^2} \bigr)\) implies that \(R\) is a neutral subspace of \(\bigl(\mathfrak K^2,\lfkip_{\mathfrak{K}^2}\bigr)\) if and only if \(\bigl\llf \{f,g\}, \{f,g\} \bigr\rlf_{\mathfrak{K}^2} = 0\) for all \(\{f,g\} \in R\). 
Consequently, the equalities in \eqref{eq-dissym} imply that a linear relation \(R\) is symmetric in \(\bigl(\mathfrak K,\kip_{\mathfrak K}\bigr)\) if and only if for all  \(\{f,g\} \in R\) we have \(\operatorname{Im}\mkern0.5mu [g,f]_{\mathfrak K} = 0\), or equivalently, \([g,f]_{\mathfrak K} \in \mathbb{R}\). In particular, a symmetric linear relation is dissipative. 

The following lemma is the analog, for linear relations in a Krein space, of \cite[Lemma~2.1]{ADW2013}. There it is stated for dissipative linear relations in a Hilbert space, and the discussion preceding it relates the result to earlier work on maximal dissipative linear relations by Krein and Langer \cite{KL}. The proof below follows the same argument. The same proof gives the corresponding statement for symmetric linear relations, so we include it in the lemma.

\begin{lemma}\label{lem-mdd}
Let \(R\) be a dissipative (respectively, symmetric) linear relation in the Krein space \(\bigl(\mathfrak K,\kip_{\mathfrak K}\bigr)\). Then \(R(0) \subseteq R^{[*]}\mkern-2.5mu(0)\). If \(R\) is maximal dissipative (respectively, maximal symmetric),  then
\(R(0) = R^{[*]}\mkern-2.5mu(0)\).
\end{lemma}

\begin{proof}
Let \(v\in R(0)\), and let \(\{f,g\}\in R\). Since \(\{0,v\}\in R\), we have \(\{f,g+\alpha v\}\in R\) for all \(\alpha\in\mathbb C\). 
By dissipativity, \(\operatorname{Im}[g+\alpha v,f]_{\mathfrak K} \geq 0\). Hence,  
\[
- \operatorname{Im}\bigl(\alpha[v,f]_{\mathfrak K}\bigr) \leq 
\operatorname{Im}[g,f]_{\mathfrak K}
 \quad \text{for all} \quad  \alpha\in\mathbb C.
\]
Since \(\alpha\in\mathbb C\) is arbitrary, it follows that \([v,f]_{\mathfrak K}=0\). 
As \(f \in \dom R\) is arbitrary, \eqref{eq-R*(0)} implies \(v\in R^{[*]}\mkern-2.5mu(0)\). This proves \(R(0)\subseteq R^{[*]}\mkern-2.5mu(0)\). This inclusion for a symmetric linear relation is trivial.

Assume now that \(R\) is maximal dissipative. Let \(v\in R^{[*]}\mkern-2.5mu(0)\), and set
\[
S = R+\lspan\bigl\{\{0,v\}\bigr\}, 
\]
where \(+\) denotes the sum of subspaces of \(\mathfrak K^2\).
Then \(R\subseteq S\). We prove that \(S\) is dissipative. Every element of \(S\) has the form \(\{f,g+\alpha v\}\), for some \(\{f,g\}\in R\) and \(\alpha\in\mathbb C\). Since \(v\in R^{[*]}\mkern-2.5mu(0)\), by \eqref{eq-R*(0)} we have \([v,f]_{\mathfrak K}=0\) for all \(f\in\dom R\). Therefore, since \(R\) is dissipative, \(\operatorname{Im}[g+\alpha v,f]_{\mathfrak K} = \operatorname{Im}[g,f]_{\mathfrak K} \geq 0\). Thus \(S\) is dissipative. As \(R\) is maximal dissipative, we have \(S=R\). Since \(\{0,v\}\in S\), it follows that \(\{0,v\}\in R\), and hence \(v\in R(0)\). Therefore \(R^{[*]}\mkern-2mu(0)\subseteq R(0)\),  proving the equality.

Replacing ``dissipative'' by ``symmetric'' and ``\(\geq\)'' by ``\(=\)'' in the preceding paragraph yields the proof of the last inclusion for a maximal symmetric linear relation.
\end{proof}

The geometric viewpoint presented in this subsection is used extensively in \cite[Chapter~2]{AzB}, although the emphasis there is primarily on operators rather than general linear relations; see, in particular, Theorem~2.2.6, Propositions~2.3.5 and~2.3.6 therein.

The following theorem is obtained by translating Theorem~\ref{th-maxsub} into the language of linear relations in the Krein space
\(\bigl(\mathfrak K,\kip_{\mathfrak K}\bigr)\), using the dictionary in the preceding table. For a linear relation \(R\) in \(\mathfrak K\), we denote by \(-R\) its opposite relation,
\[
-R = \Bigl\{\{f,-g\}\in\mathfrak K^2:\{f,g\}\in R\Bigr\}.
\]
Clearly \(R\) is closed if and only if \(-R\) is closed. By the definition of the Lagrange indefinite inner product in \eqref{bdHf}, the subspace \(R\) of \(\bigl(\mathfrak K^2,\lfkip_{\mathfrak K^2}\bigr)\) is nonpositive if and only if \(-R\) is nonnegative. Equivalently, \(-R\) is dissipative in \(\bigl(\mathfrak K,\kip_{\mathfrak K}\bigr)\) if and only if \(R\), viewed as a subspace of \(\bigl(\mathfrak K^2,\lfkip_{\mathfrak K^2}\bigr)\), is nonpositive. Moreover, \((-R)^{[*]}=-R^{[*]}\). Combining these observations with Theorem~\ref{th-maxsub}, we obtain the following theorem, see also \cite[Lemma~1.4]{ACD16}.

\begin{theorem} \label{th-maxnnR}
Let \((\mathfrak{K},\kip_\mathfrak{K})\) be a Krein space and let \(R\) be a linear relation in \(\mathfrak{K}\). The following statements hold.
\begin{enumerate}
\renewcommand*\theenumi{\alph{enumi}}
\renewcommand*\labelenumi{\rm (\theenumi)}
\item  \label{th-maxnnR-i1}
\(R\) is maximal dissipative if and only if \(R\) is closed and both \(R\) and \(-R^{[*]}\) are dissipative.
\item  \label{th-maxnnR-i2}
\(R\) is maximal symmetric if and only if \(R\) is symmetric and \(R\) or \(-R\) is maximal dissipative.
\item  \label{th-maxnnR-i3}
\(R\) is self-adjoint if and only if \(R\) and \(-R\) are maximal dissipative.
\end{enumerate}
\end{theorem}

\section{Compressions in Krein spaces} \label{sec-comps}

Let \(\bigl(\widetilde{\mathfrak H},\kip_{\widetilde{\mathfrak H}}\bigr)\) be a Krein space. Let \(\mathfrak{H}\) and \(\widehat{\mathfrak{H}}\) be Krein subspaces of \(\widetilde{\mathfrak H}\) which are orthogonal complements of each other. Then
\begin{equation}\label{eq-HhHds}
\widetilde{\mathfrak H} = \mathfrak{H} [\dot{+}] \widehat{\mathfrak{H}},
\end{equation}
where the sum is direct and orthogonal in \(\kip_{\widetilde{\mathfrak H}}\).
We denote by \(\widetilde{P}_{\mathfrak H}\) and \(\widetilde{P}_{\widehat{\mathfrak H}}\) the corresponding orthogonal projections. When \(\mathfrak H\) and \(\widehat{\mathfrak{H}}\) are considered as Krein spaces in their own right, their inner products will be denoted by \(\kip_{\mathfrak H}\) and \(\kip_{\widehat{\mathfrak{H}}}\), respectively; these are the restrictions of \(\kip_{\widetilde{\mathfrak H}}\) to the corresponding subspaces.
Let \(J\) be a fundamental symmetry on \(\mathfrak H\) and let \(\widehat J\) be a  fundamental symmetry on \(\widehat{\mathfrak{H}}\). Let \(\ahip_{\mathfrak H}\) and \(\ahip_{\widehat{\mathfrak{H}}}\) be the associated Hilbert space inner products, and denote the corresponding norms by \(\nor_{\mathfrak H}\) and \(\nor_{\widehat{\mathfrak{H}}}\). Then \(\widetilde J = J\oplus\widehat J\) is a fundamental symmetry on \(\widetilde{\mathfrak H}\). We denote the corresponding Hilbert space inner product by \(\ahip_{\widetilde{\mathfrak H}}\) and by \(\nor_{\widetilde{\mathfrak{H}}}\) its norm. Also with respect to this Hilbert space inner product, the decomposition \eqref{eq-HhHds} is orthogonal.

Let \(R\) be a linear relation in \(\widetilde{\mathfrak H}\), that is, a subspace \(R \subseteq \widetilde{\mathfrak{H}}^2\). The \emph{restriction of} \(R\) \emph{to} \(\mathfrak{H}\), denoted by \(R|_{\mathfrak H}\), is the linear relation with the domain in \(\mathfrak{H}\) and the codomain in \(\widetilde{\mathfrak H}\) defined by
\[
R|_{\mathfrak H}
= R \cap \bigl(\mathfrak{H}\times\widetilde{\mathfrak H}\bigr)
= \Bigl\{ \{f, \widetilde{g}\} \in R : f \in \mathfrak{H}\Bigr\}.
\]
The relation \(\widetilde{P}_{\mathfrak H} R = \bigl\{\{\widetilde{f}, \widetilde{P}_{\mathfrak H} \widetilde{g}\}:\{\widetilde{f},\widetilde{g}\}\in R \bigr\}\) is the \emph{composition} of relations \(\widetilde{P}_{\mathfrak H}\) and \(R\). The \emph{compression of} \(R\) \emph{to} \(\mathfrak{H}\), denoted by \(\CH(R)\), is the linear relation in \(\mathfrak H\) defined by 
\[
\CH(R) = \widetilde{P}_{\mathfrak H} \, R|_{\mathfrak H}
= \Bigl\{ \bigl\{f, \widetilde{P}_{\mathfrak H} \widetilde{g} \bigr\} : \{f, \widetilde{g}\} \in R|_{\mathfrak H} \Bigr\}.
\]

We apply the construction in Subsection~\ref{subs-ddip}
to \(\widetilde{\mathfrak{H}}\) with the fundamental symmetry
\(\widetilde{J}\).
On \(\widetilde{\mathfrak{H}}^2\), we denote the resulting
Lagrange inner product by \(\lfkip\), the fundamental symmetry
by \(\widetilde{\mathbb{J}}\), and the corresponding Hilbert
space inner product by \(\dahip\).
The Hilbert space \(\bigl(\widetilde{\mathfrak{H}}^2,\dahip\bigr)\) is the direct product of two copies of the Hilbert space \(\bigl(\widetilde{\mathfrak{H}},\ahip_{\widetilde{\mathfrak{H}}}\bigr)\).

Notice that \(\mathfrak H^2\) and \(\widehat{\mathfrak H}^{2}\) are Krein subspaces of \(\bigl(\widetilde{\mathfrak{H}}^2, \lfkip\bigr)\). The Krein subspace \(\bigl(\mathfrak{H}^2, \lfkip\bigr)\) is the Lagrange-Krein space associated with the Krein space \(\bigl(\mathfrak{H},\kip_{\mathfrak{H}}\bigr)\).  The subspace \(\mathfrak H\times\widetilde{\mathfrak H}\), which is used in the definition of the restriction \(R|_{\mathfrak H}\), is closed but degenerate when \(\widehat{\mathfrak H}\ne\{0\}\); its isotropic part is \(\{0\}\times\widehat{\mathfrak H}\). Similarly, the subspace \(\widetilde{\mathfrak{H}}\times\mathfrak{H}\), which plays a prominent role in Section~\ref{SecCAACn}, is also closed, and it is degenerate when \(\widehat{\mathfrak H}\ne\{0\}\); its isotropic part is \(\widehat{\mathfrak H}\times\{0\}\).

We denote adjoints in \(\bigl(\widetilde{\mathfrak{H}},\kip_{\widetilde{\mathfrak{H}}}\bigr)\)
and \(\bigl(\mathfrak{H},\kip_{\mathfrak{H}}\bigr)\) by \(^{\widetilde{*}}\) and \(^{*}\), respectively. 
Recall that the adjoint \(R^{\widetilde{*}}\) of a linear relation \(R\) in \(\bigl(\widetilde{\mathfrak H},\kip_{\widetilde{\mathfrak H}}\bigr)\) is its orthogonal complement
\(R^{\llf\perp\rlf}\) in \(\bigl(\widetilde{\mathfrak{H}}^2,\lfkip\bigr)\).
In particular, a linear relation \(T\) in \(\bigl(\mathfrak{H},\kip_{\mathfrak{H}}\bigr)\) can also be regarded as a linear relation in \(\bigl(\widetilde{\mathfrak{H}},\kip_{\widetilde{\mathfrak{H}}}\bigr)\), and its two adjoints are related by 
\(T^{\widetilde{*}}=T^*\llf\dotplus\rlf\widehat{\mathfrak{H}}^2\).

\begin{proposition} \label{pr-inher}
Let \(R\) be a linear relation in the Krein space \(\widetilde{\mathfrak H} = \mathfrak{H} [\dot{+}] \widehat{\mathfrak H}\). The following statements hold.
\begin{enumerate}
\renewcommand*\theenumi{\roman{enumi}}
\renewcommand*\labelenumi{{\rm (\theenumi)}}
\item \label{pr-inher-i1}
The inclusion \(\CH(R^{\widetilde{*}})\subseteq \CH(R)^*\) holds.

\item  \label{pr-inher-i2}
If \(R\) is dissipative (respectively, symmetric) in \(\bigl(\widetilde{\mathfrak{H}}, \kip_{\widetilde{\mathfrak{H}}} \bigr)\), then \(\CH(R)\) is dissipative (respectively, symmetric) in \(\bigl(\mathfrak{H}, \kip_{\mathfrak{H}} \bigr)\).
\end{enumerate}
\end{proposition}

\begin{proof}

(\ref{pr-inher-i1})
By the definitions of the adjoint and the compression, for all \(\{u, v\} \in \mathfrak{H}^2\) we have  
\begin{equation}
\label{eq-CHR*}
  \{u,v\} \in \CH(R)^* \quad \Leftrightarrow \quad
\bllf \{f,\widetilde{g}\}, \{u,v\}\brlf = 0 \quad \forall \{f,\widetilde{g}\} \in R|_{\mathfrak H}.
\end{equation}
Let \(\{u,v\}\in\CH(R^{\widetilde{*}})\) be arbitrary.
Then there exists \(\widehat{v}\in\widehat{\mathfrak{H}}\)
such that \(\{u,v+\widehat{v}\}\in R^{\widetilde{*}}\).
Since \(R|_{\mathfrak H}\subseteq R\), the orthogonality
of \(\mathfrak H\) and \(\widehat{\mathfrak H}\) gives
\[
\bllf\{f,\widetilde{g}\},\{u,v\}\brlf
=
\bllf\{f,\widetilde{g}\},\{u,v+\widehat{v}\}\brlf
=0
\quad\text{for all}\quad
\{f,\widetilde{g}\}\in R|_{\mathfrak H}.
\]
By \eqref{eq-CHR*}, we have \(\{u,v\}\in\CH(R)^*\).
This proves (\ref{pr-inher-i1}).

(\ref{pr-inher-i2})
Let \(\{f,h\}\in\CH(R)\), and choose
\(\{f,\widetilde{g}\}\in R|_{\mathfrak H}\) such that
\(h=\widetilde{P}_{\mathfrak H}\widetilde{g}\).
Since \(f\in\mathfrak H\), the orthogonality of \(\mathfrak H\) and
\(\widehat{\mathfrak H}\) gives
\[
\bllf\{f,h\},\{f,h\}\brlf
=
\bllf\{f,\widetilde{g}\},\{f,\widetilde{g}\}\brlf.
\]
If \(R\) is dissipative (respectively, symmetric), the
right-hand side is nonnegative (respectively, zero).
Hence \(\CH(R)\) is dissipative (respectively, symmetric).
\end{proof}

\begin{example}  \label{ex-com-op}
Let \(R\) be a linear relation in \(\widetilde{\mathfrak{H}}\).  The orthogonal decomposition
\(\widetilde{\mathfrak{H}}
=\mathfrak{H}[\dotplus]\widehat{\mathfrak{H}}\) is said to \emph{reduce} \(R\) if
\[
R=\bigl(R\cap\mathfrak{H}^2\bigr)
\llf\dotplus\rlf
\bigl(R\cap\widehat{\mathfrak{H}}^2\bigr).
\]
Clearly, if \(R\) is reduced by \(\widetilde{\mathfrak{H}}
=\mathfrak{H}[\dotplus]\widehat{\mathfrak{H}}\), then \(\CH(R)=R\cap\mathfrak{H}^2\).  

A special case of reduction occurs when \(R\) is a closed linear relation in
\(\widetilde{\mathfrak{H}}\) such that
\(R(0)=R^{\widetilde{*}}\mkern-2mu(0)=\widehat{\mathfrak{H}}\).
Indeed, by \eqref{eq-R*(0)} and \(R^{\widetilde{*}}(0)=\widehat{\mathfrak{H}}\)
we have
\(\overline{\operatorname{dom}R}=\mathfrak{H}\).
In particular, \(\operatorname{dom}R\subseteq\mathfrak{H}\). 
For every \(\{f,\widetilde{g}\}\in R|_{\mathfrak{H}}\), the equality
\(R(0)=\widehat{\mathfrak{H}}\) implies that
\(\{0,\widetilde{P}_{\widehat{\mathfrak{H}}}\widetilde{g}\}\in R\).
Subtracting this pair gives
\(\{f,\widetilde{P}_{\mathfrak{H}}\widetilde{g}\}\in R\cap\mathfrak{H}^2\).
Consequently,
\begin{equation} \label{eq-need}  
\CH(R)=R\cap\mathfrak{H}^2  \quad \text{and} 
\quad
R=\CH(R)\dotplus
\bigl(\{0\}\times\widehat{\mathfrak{H}}\bigr).
\end{equation} 
The linear relation \(\CH(R)\) is closed, its domain is
\(\operatorname{dom}R\), and its multivalued part is
\(R(0)\cap\mathfrak{H}=\{0\}\).
Thus \(\CH(R)\) is a closed, densely defined operator in
\(\mathfrak{H}\), and \eqref{eq-need} shows that it is the operator part of \(R\). 
\end{example}

\begin{example} 
\label{ex-compression-matrix}
As in Subsection~\ref{subs-ddip}, let \(\bigl(\mathfrak{K},\kip_\mathfrak{K} \bigr)\)  be a Krein space with a fundamental symmetry \(J\) on it. Let \(\widetilde{\mathfrak{H}} = \mathfrak{K}^2\) be equipped with the direct product Krein space inner product. That is,  
\begin{equation*}  
\bigl[ \{f,g\}, \{h,k\} \bigr]_{\mathfrak{K}^2} = 
[ f, h ]_\mathfrak{K} + [g,k]_\mathfrak{K}, \qquad \{f,g\}, \{h,k\} \in \mathfrak{K}^2.
\end{equation*}
Then \(\bigl(\mathfrak{K}^2,\kip_{\mathfrak{K}^2}\bigr) = \bigl(\widetilde{\mathfrak{H}}, \kip_{\widetilde{\mathfrak{H}}}\bigr)\) is a Krein space and \(\mathfrak{H} = \mathfrak{K} \times \{0\}\) and \(\widehat{\mathfrak{H}} = \{0\}\times\mathfrak{K}\) are Krein subspaces of \(\widetilde{\mathfrak{H}} = \mathfrak{K}^2\). 

Let \(\mathbb{J}\) be the everywhere defined linear involution defined in \eqref{eq-dbbJ}. Here we represent \(\mathbb{J}\) as its graph: 
\[
\mathbb{J} = \Bigl\{ \bigl\{ \{f,g\}, \{-\iu J g, \iu J f\}  \bigr\}  :  
\{f,g\} \in  \mathfrak{K}^2 \Bigr\}. 
\]
Then \(\mathbb{J}\) is self-adjoint in \(\bigl(\widetilde{\mathfrak{H}}, \kip_{\widetilde{\mathfrak{H}}}\bigr)\) and 
\[
\mathbb{J} \cap \mathfrak{H}^2 =  \Bigl\{  \bigl\{\{0,0\},  \{0,0\} \bigr\} \Bigr\}, 
\qquad 
\mathbb{J}|_{\mathfrak{H}} = \Bigl\{ 
\bigl\{ \{f, 0\}, \{ 0 , \iu J f\}  \bigr\}  :  
f \in  \mathfrak{K} \Bigr\}, 
\]
and 
\[
\CH(\mathbb{J}) = \Bigl\{ 
\bigl\{ \{f, 0\}, \{ 0 , 0\}  \bigr\}  :  
f \in  \mathfrak{K} \Bigr\} = \mathfrak{H} \times \bigl\{ \{0,0\}\bigr\}. 
\]
That is, \(\CH(\mathbb{J})\) is the zero operator defined on \(\mathfrak{H}\), while \(\mathbb{J} \cap \mathfrak{H}^2\) consists of only the zero vector in \(\mathfrak{H}^2\). 
\end{example}

\begin{example} 
\label{ex-compression-multivalued}
Let \(\bigl(\mathfrak{K},\kip_\mathfrak{K} \bigr)\) be a Krein space. Let \(\widetilde{\mathfrak{H}} = \mathfrak{K}^2\) be equipped with the direct product Krein space inner product of the Krein spaces \(\bigl(\mathfrak{K},\kip_\mathfrak{K} \bigr)\) and \(\bigl(\mathfrak{K},-\kip_\mathfrak{K} \bigr)\). That is,  
\begin{equation*}  
\bigl[ \{f,g\}, \{h,k\} \bigr]_{\mathfrak{K}^2} = 
[ f, h ]_\mathfrak{K} - [g,k]_\mathfrak{K}, \qquad \{f,g\}, \{h,k\} \in \mathfrak{K}^2.
\end{equation*}
Then \(\bigl(\mathfrak{K}^2,\kip_{\mathfrak{K}^2}\bigr) = \bigl(\widetilde{\mathfrak{H}}, \kip_{\widetilde{\mathfrak{H}}}\bigr)\) is a Krein space and \(\mathfrak{H} = \mathfrak{K} \times \{0\}\) and \(\widehat{\mathfrak{H}} = \{0\}\times\mathfrak{K}\) are Krein subspaces of \(\widetilde{\mathfrak{H}} = \mathfrak{K}^2\). 

Define a linear relation \(R\) on \(\widetilde{\mathfrak{H}} = \mathfrak{K}^2\) by 
\[
R = \Bigl\{ \bigl\{\widetilde{f}, \widetilde{g} \bigr\} \in \widetilde{\mathfrak{H}}^2  : \widetilde{f} = \{f,f\}, \  \widetilde{g}=\{g,g\}, \ f,g  \in  \mathfrak{K} \Bigr\}.
\]
We show that \(R\) is self-adjoint in \(\widetilde{\mathfrak{H}}\).
Let \(\widetilde{x}=\{x_1,x_2\}\) and
\(\widetilde{y}=\{y_1,y_2\}\) be arbitrary elements of
\(\widetilde{\mathfrak{H}}\). Then
\begin{align*}
\bigl\{\widetilde{x},\widetilde{y}\bigr\}
\in R^{\widetilde{*}}
&\quad\Leftrightarrow\quad
[f,y_1-y_2]_{\mathfrak{K}}
-[g,x_1-x_2]_{\mathfrak{K}}=0
\quad\text{for all}\quad f,g\in\mathfrak{K}
\\
&\quad\Leftrightarrow\quad
x_1=x_2,\quad y_1=y_2
\\
&\quad\Leftrightarrow\quad
\bigl\{\widetilde{x},\widetilde{y}\bigr\}\in R.
\end{align*}
The second equivalence follows from nondegeneracy of \(\kip_{\mathfrak{K}}\).

Furthermore,
\[
R\cap\mathfrak{H}^2
=
\Bigl\{
\bigl\{\{0,0\},\{0,0\}\bigr\}
\Bigr\}, 
\qquad
R|_{\mathfrak{H}}
=
\Bigl\{
\bigl\{\{0,0\},\{g,g\}\bigr\}:
g\in\mathfrak{K}
\Bigr\},
\]
and 
\[
\CH(R) = \Bigl\{ 
\bigl\{ \{0, 0\}, \{ g , 0\}  \bigr\}  :  
g \in  \mathfrak{K} \Bigr\} =  \bigl\{ \{0,0\} \bigr\} \times \mathfrak{H}. 
\]
Thus \(\CH(R)\) is a purely multivalued relation:
its domain is \(\{0\}\) and its multivalued part is
\(\mathfrak{H}\).
\end{example}

\section{Compression commutes with taking adjoints} \label{SecCAACn}

In this section we give a sufficient condition under which a closed
linear relation \(R\) in \(\widetilde{\mathfrak{H}}\) satisfies
\(\CH(R^{\widetilde{*}})=\CH(R)^*\), that is, its compression
commutes with taking adjoints.   
Our starting point is the part in the proof of \cite[Lemma~4.4]{ADW2013} where it is shown that if \(\widetilde{\mathfrak{H}}\) is a Hilbert space and the linear relations \(R\) and \(\widetilde{P}_{\mathfrak{H}}R^{\widetilde{*}}\) are closed, then the adjoint identity holds. We establish this implication for arbitrary closed linear relations \(R\) in Krein spaces. We show that the closedness of \(\widetilde{P}_{\mathfrak{H}}R^{\widetilde{*}}\) is equivalent to  the closedness of the subspace \(\widetilde{P}_{\widehat{\mathfrak{H}}}(\dom R)\). This implies that when the exit space \(\widehat{\mathfrak{H}}\) is finite dimensional, both closedness conditions and the adjoint identity hold automatically.

\begin{lemma} \label{le-prs}
Let \(\bigl( \mathfrak{G}, \ahip \bigr) \) be a Hilbert space and let \(\mathfrak{A}\) and \(\mathfrak{B}\) be closed subspaces of \(\mathfrak{G}\). Then \(P_{\mathfrak{A}}\mathfrak{B}\) is closed if and only if \(P_{\mathfrak{A}^{\mkern-2mu\langle\mkern-2mu\perp\mkern-2mu\rangle}}%
\mathfrak{B}^{\mkern-1mu\langle\mkern-1mu\perp\mkern-1mu\rangle}\) is closed. Here \(P_{\mathfrak{A}}\),  \(P_{\mathfrak{A}^{\mkern-2mu\langle\mkern-1mu\perp\mkern-1mu\rangle}}\) are the orthogonal projections onto the indicated subspaces.
\end{lemma}

\begin{proof}
Since \(\mathfrak{A}\) is closed, and
\begin{equation*}
\mathfrak{A} + \mathfrak{B}^{\mkern-1mu\langle\mkern-1mu\perp\mkern-1mu\rangle} =
\mathfrak{A} \langle\mkern-1mu\dotplus\mkern-1mu\rangle P_{\mathfrak{A}^{\mkern-2mu\langle\mkern-1mu\perp\mkern-1mu\rangle}} \mathfrak{B}^{\mkern-1mu\langle\mkern-1mu\perp\mkern-1mu\rangle},
\end{equation*}
the subspace \(\mathfrak{A} + \mathfrak{B}^{\mkern-1mu\langle\mkern-1mu\perp\mkern-1mu\rangle}\) 
is closed if and only if 
\(P_{\mathfrak{A}^{\mkern-2mu\langle\mkern-1mu\perp\mkern-1mu\rangle}} \mathfrak{B}^{\mkern-1mu\langle\mkern-1mu\perp\mkern-1mu\rangle}\) 
is closed. Since 
\(\mathfrak{A}^{\mkern-1mu\langle\mkern-1mu\perp\mkern-1mu\rangle}\) 
is closed, and
\begin{equation*}
\mathfrak{A}^{\mkern-1mu\langle\mkern-1mu\perp\mkern-1mu\rangle} + \mathfrak{B} =
\mathfrak{A}^{\mkern-1mu\langle\mkern-1mu\perp\mkern-1mu\rangle} \langle\mkern-1mu\dotplus\mkern-1mu\rangle P_{\mathfrak{A}} \mathfrak{B},
\end{equation*}
the subspace \(\mathfrak{A}^{\mkern-1mu\langle\mkern-1mu\perp\mkern-1mu\rangle} + \mathfrak{B}\) is closed if and only if \(P_{\mathfrak{A}} \mathfrak{B}\)  is closed.

By \cite[Theorem~IV.4.8]{Kato} the sum of two closed subspaces is closed if and only if the sum of their orthogonal complements is closed. Therefore \(\mathfrak{A} + \mathfrak{B}^{\mkern-1mu\langle\mkern-1mu\perp\mkern-1mu\rangle}\) is closed if and only if \(\mathfrak{A}^{\mkern-1mu\langle\mkern-1mu\perp\mkern-1mu\rangle} + \mathfrak{B}\) is closed. These three equivalences prove the lemma.
\end{proof}

\begin{theorem}\label{th-domRandPR*}
Let \(R\) be a closed linear relation in the Krein space \(\widetilde{\mathfrak H} = \mathfrak{H}[\dot{+}] \widehat{\mathfrak H} \). The linear relation
\(\widetilde{P}_{\mathfrak H}R^{\widetilde{*}}\)
is a closed linear relation in \(\widetilde{\mathfrak H}\) if and only if the subspace
\(\widetilde{P}_{\widehat{\mathfrak H}}(\dom R)\) is closed in \(\widehat{\mathfrak H}\).
\end{theorem}

\begin{proof}
We use the notation introduced in Section~\ref{sec-comps}. In particular we use the Hilbert space inner product \(\dahip\) on \(\widetilde{\mathfrak{H}}^2\) and the fact that \(\widetilde{P}_{\widehat{\mathfrak H}}\widetilde{J} = \widetilde{J} \widetilde{P}_{\widehat{\mathfrak H}}\).

We apply Lemma~\ref{le-prs} to the Hilbert space \(\bigl(\mathfrak{G}, \ahip \bigr) =  \bigl(\widetilde{\mathfrak{H}}^2,\dahip\bigr)\) and its closed subspaces \(\mathfrak{A} = \widetilde{\mathfrak H}\times\mathfrak{H}\) and \(\mathfrak{B} = R^{\widetilde{*}}\).
Then \(\mathfrak{A}^{\mkern-1mu\langle\mkern-1mu\perp\mkern-1mu\rangle} 
= \{0\}\times \widehat{\mathfrak{H}}\) and \(P_{\mathfrak{A}}\mathfrak{B} = \widetilde{P}_{\mathfrak{H}} R^{\widetilde{*}}\).

Since the environment of this proof is the Hilbert space \(\bigl(\widetilde{\mathfrak{H}}^2,\dahip\bigr)\), it is convenient to write \(R^{\widetilde{*}}\) as \(R^{\widetilde{*}} = R^{\llf\perp\rlf} = \bigl(\widetilde{\mathbb{J}}R\bigr)^{\mkern-2mu\lda\mkern-1mu\perp\mkern-1mu\rda}\). 
Then \(\ran \bigl(\widetilde{\mathbb{J}}R\bigr) = \widetilde{J}(\dom R)\), compare \eqref{eq-dbbJ}; \(\widetilde{\mathbb{J}}R \) is closed since \(R\) is closed; and  
\((R^{\widetilde{*}})^{\mkern-1mu\lda\mkern-1mu\perp\mkern-1mu\rda} 
= \widetilde{\mathbb{J}}R\). Further,
\begin{align*}
P_{\mathfrak{A}^{\mkern-2mu\langle\mkern-1mu\perp\mkern-1mu\rangle}} \mathfrak{B}^{\mkern-1mu\langle\mkern-1mu\perp\mkern-1mu\rangle}
& =
\Bigl\{ \bigl\{0,\widetilde{P}_{\widehat{\mathfrak{H}}} \widetilde{y}\bigr\} : 
\bigl\{ \widetilde{x}, \widetilde{y} \bigr\} \in (R^{\widetilde{*}})^{\lda\perp\rda}\Bigr\} \\
& = \{0\} \times \widetilde{P}_{\widehat{\mathfrak H}}%
\bigl(\ran\bigl(\widetilde{\mathbb{J}}R\bigr)\bigr) \\
& = \{0\} \times  \widetilde P_{\widehat{\mathfrak H}} \widetilde{J}(\dom R) \\
& = \{0\} \times   \widetilde{J}\widetilde P_{\widehat{\mathfrak H}}(\dom R).
\end{align*}

By Lemma~\ref{le-prs} we have that \(P_{\mathfrak{A}}\mathfrak{B} = \widetilde{P}_{\mathfrak{H}} R^{\widetilde{*}}\) is closed if and only if \(P_{\mathfrak{A}^{\mkern-1mu\langle\mkern-1mu\perp\mkern-1mu\rangle}}% 
\mathfrak{B}^{\mkern-1mu\langle\mkern-1mu\perp\mkern-1mu\rangle} 
= \{0\} \times \widetilde{J}\widetilde P_{\widehat{\mathfrak H}}(\dom R)\) is closed. 
Clearly \(\{0\} \times \widetilde{J}\widetilde P_{\widehat{\mathfrak H}}(\dom R)\) is closed if and only if \(\widetilde{J}\widetilde P_{\widehat{\mathfrak H}}(\dom R)\) is closed. Finally, since \(\widetilde{J}\) is a bounded involution, the last space is closed if and only if \(\widetilde P_{\widehat{\mathfrak H}}(\dom R)\) is closed. This sequence of equivalences proves the theorem.
\end{proof}

Since for every linear relation \(R\) in \(\widetilde{\mathfrak{H}}\) we have \(\widebar{R} = (R^{\widetilde{*}})^{\widetilde{*}}\), applying the preceding theorem to the closed linear relation \(R^{\widetilde{*}}\) yields the following statement. 

\begin{corollary} \label{co-domRandPR*c1}
Let \(R\) be a linear relation in the Krein space \(\widetilde{\mathfrak H} = \mathfrak{H}[\dot{+}] \widehat{\mathfrak H} \). The relation
\(\widetilde{P}_{\mathfrak H} \widebar{R}\) is closed if and only if the subspace
\(\widetilde{P}_{\widehat{\mathfrak H}}\bigl(\dom R^{\widetilde{*}}\bigr)\) is closed in \(\widehat{\mathfrak H}\).
\end{corollary}

\begin{theorem} \label{th-comadjn}
Let \(\widetilde{\mathfrak H} = \mathfrak{H} [\dot{+}] \widehat{\mathfrak H}\) be a Krein space and let \(R\) be a closed linear relation in  \(\widetilde{\mathfrak{H}}\).  If \(\widetilde{P}_{\mathfrak H}R^{\widetilde{*}}\) is a closed linear relation, or equivalently if \(\widetilde{P}_{\widehat{\mathfrak H}}(\dom R)\) is closed in \(\widehat{\mathfrak H}\), then \(\CH(R^{\widetilde{*}})\) is closed and \(\CH(R^{\widetilde{*}}) = \CH(R)^*\). 
\end{theorem}

\begin{proof}
The proof is an adaptation of the product-adjoint argument used in the proof of  \cite[Lemma~4.4]{ADW2013}. Put \(P = \widetilde P_{\mathfrak H}\). Thus \(P\) is a bounded self-adjoint projection in the Krein space \(\widetilde{\mathfrak H}\). We use the following product rule for adjoints of linear relations, which is a special case of \cite[Lemma~2.9(i)]{DHMS2009}:
if \(A\) is a bounded, everywhere defined operator and \(B\) is a linear relation in a Krein space \(\bigl(\mathfrak{K},\kip\bigr)\), then \((AB)^{[*]}=B^{[*]}A^{[*]}\).

Assume that \(R\) and \(PR^{\widetilde{*}}\) are closed linear relations in \(\widetilde{\mathfrak H}\). Since \(R\) is closed, applying the product rule gives \(RP = (R^{\widetilde{*}})^{\widetilde{*}}P^{\widetilde{*}} = (PR^{\widetilde{*}})^{\widetilde{*}}\). Using this, the closedness of \(PR^{\widetilde{*}}\) and the product rule once more we get
\begin{equation}\label{eq-PRP}
(PRP)^{\widetilde{*}}
 =
\bigl(P(PR^{\widetilde{*}})^{\widetilde{*}}\bigr)^{\widetilde{*}}
 =
(PR^{\widetilde{*}})^{\widetilde{*}\widetilde{*}}P^{\widetilde{*}}
=
(PR^{\widetilde{*}})P
=
PR^{\widetilde{*}}P.
\end{equation}

From the definition of \(PRP\) it follows that \(PRP = \CH(R)\mathbin{\dot+}(\widehat{\mathfrak{H}}\times\{0\})\), where \(\dotplus\) denotes a direct sum of subspaces of \(\widetilde{\mathfrak{H}}^2\). The decomposition 
\(\widetilde{\mathfrak{H}}^2 
=\mathfrak{H}^2\llf\dotplus\rlf \widehat{\mathfrak{H}}^2\)
is orthogonal with respect to \(\lfkip\), and \(\widehat{\mathfrak{H}}\times\{0\}\) is hyper-maximal neutral in \(\bigl(\widehat{\mathfrak{H}}^2,\lfkip\bigr)\). 
Taking orthogonal complements in \(PRP = \CH(R)\mathbin{\dot+}(\widehat{\mathfrak{H}}\times\{0\})\) therefore gives
\[
(PRP)^{\widetilde{*}}
=\CH(R)^*\mathbin{\dotplus}
\bigl(\widehat{\mathfrak{H}}\times\{0\}\bigr).
\] 

Since also \((PR^{\widetilde{*}}P) \cap \mathfrak{H}^2=\CH(R^{\widetilde{*}})\), the equality in \eqref{eq-PRP} proves \(\CH(R^{\widetilde{*}}) = \CH(R)^*\).

As \(\CH(R)^*\) is an adjoint linear relation in the Krein space \(\mathfrak H\), it is closed. Therefore \(\CH(R^{\widetilde{*}})\) is closed as well.
\end{proof}

The closedness condition in Theorem~\ref{th-comadjn} is sufficient but not necessary for compression to commute with taking adjoints, as the following example shows.

\begin{example} \label{exa-twocon}
Let \(\widetilde{\mathfrak{H}}
=\mathfrak{H}[\dotplus]\widehat{\mathfrak{H}}\) be a Hilbert
space such that \(\mathfrak{H}\neq\{0\}\) and
\(\widehat{\mathfrak{H}}\) is infinite dimensional.
Let \(\widehat{A}\) be an unbounded self-adjoint operator
in \(\widehat{\mathfrak{H}}\), and let
\(R=0_{\mathfrak{H}}\dotplus\widehat{A}\) be the direct sum
of the zero operator \(0_{\mathfrak{H}}\) and \(\widehat{A}\).
Then \(R=R^{\widetilde{*}}\), and the orthogonal decomposition
of \(\widetilde{\mathfrak{H}}\) reduces \(R\).
By the first part of Example~\ref{ex-com-op}, we have
\(\CH(R^{\widetilde{*}})=\CH(R)^*=0_{\mathfrak{H}}\).

On the other hand,
\(\widetilde{P}_{\mathfrak{H}}R^{\widetilde{*}}\) is the
zero operator with domain
\(\mathfrak{H}[\dotplus](\dom\widehat{A})\).
Since \(\widehat{A}\) is unbounded and self-adjoint,
its domain is a dense proper subspace of
\(\widehat{\mathfrak{H}}\), and hence is not closed.
Consequently,
\(\widetilde{P}_{\mathfrak{H}}R^{\widetilde{*}}\)
is not closed. Equivalently, the projected-domain condition
fails because
\(\widetilde{P}_{\widehat{\mathfrak{H}}}(\dom R)
=\dom\widehat{A}\) is not closed in
\(\widehat{\mathfrak{H}}\).
The setting demonstrated in this example is further clarified in  Appendix~\ref{appendix}. 
\end{example}

Applying Theorem~\ref{th-comadjn} to \(R^{\widetilde{*}}\) yields the following corollary.

\begin{corollary} \label{co-thR*}
Let \(\widetilde{\mathfrak H} = \mathfrak{H} [\dot{+}] \widehat{\mathfrak H}\) be a Krein space and let \(R\) be a linear relation in  \(\widetilde{\mathfrak{H}}\). If \(\widetilde{P}_{\widehat{\mathfrak H}}\bigl(\dom R^{\widetilde{*}}\bigr)\) is closed in \(\widehat{\mathfrak H}\),  or equivalently if \(\widetilde{P}_{\mathfrak H}\widebar{R}\) is a closed linear relation in \(\widetilde{\mathfrak H}\), then \(\CH(\widebar{R})\) is closed and \(\CH(\widebar{R}) = \CH(R^{\widetilde{*}})^*\).
\end{corollary}

\section{Properties preserved and reflected under compression} \label{SecPR}

We first establish preservation results, namely, results asserting that certain properties of a linear relation \(R\) are preserved by its compression \(\CH(R)\). The hypothesis in Theorem~\ref{th-comadj4} is the closedness of 
\(\widetilde P_{\widehat{\mathfrak H}}\bigl(\dom R^{\widetilde{*}}\bigr)\) in \(\widehat{\mathfrak H}\), while in Theorem~\ref{th-conv} the main hypothesis is the closedness of \(\widetilde P_{\widehat{\mathfrak H}}(\dom R)\) in \(\widehat{\mathfrak H}\). The two projected domain conditions are independent, even for closed, densely defined operators in Hilbert spaces as the following example shows.

\begin{example}

Let \(\mathfrak{K}\) be an infinite-dimensional Hilbert
space, and equip \(\widetilde{\mathfrak{H}}=\mathfrak{K}^2\)
with the product Hilbert space inner product. Set
\(\mathfrak{H}=\mathfrak{K}\times\{0\}\) and
\(\widehat{\mathfrak{H}}=\{0\}\times\mathfrak{K}\).
Then
\(\widetilde{\mathfrak{H}}
=\mathfrak{H}[\dotplus]\widehat{\mathfrak{H}}\).

Let \(A\) be an unbounded self-adjoint operator in
\(\mathfrak{K}\), and define an operator \(R\) in
\(\widetilde{\mathfrak{H}}\) by
\[
\dom R=\mathfrak{K}\times\dom A,
\qquad
R\{f,g\}=\{Ag,0\}
\quad\text{for all}\quad
\{f,g\}\in\dom R.
\]
Since \(A\) is closed and densely defined, so is \(R\).
The self-adjointness of \(A\) gives
\[
\dom R^{\widetilde{*}}=\dom A\times\mathfrak{K},
\qquad
R^{\widetilde{*}}\{u,v\}=\{0,Au\}
\quad\text{for all}\quad
\{u,v\}\in\dom R^{\widetilde{*}}.
\]
Consequently,
\[
\widetilde{P}_{\widehat{\mathfrak{H}}}(\dom R)
=\{0\}\times\dom A,
\qquad
\widetilde{P}_{\widehat{\mathfrak{H}}}
\bigl(\dom R^{\widetilde{*}}\bigr)
=\widehat{\mathfrak{H}}.
\]
Since \(A\) is unbounded and self-adjoint, its domain is
a dense proper subspace of \(\mathfrak{K}\), and hence
is not closed. Thus the first projected domain is not
closed, whereas the second is closed.

Replacing \(R\) by \(R^{\widetilde{*}}\) reverses these
conclusions, because
\((R^{\widetilde{*}})^{\widetilde{*}}=R\).
\end{example}

\begin{theorem} \label{th-comadj4}
Let \(\widetilde{\mathfrak H} = \mathfrak{H} [\dot{+}] \widehat{\mathfrak H}\) be a Krein space and let \(R\) be a linear relation in  \(\widetilde{\mathfrak{H}}\). Assume \(\widetilde{P}_{\widehat{\mathfrak H}}\bigl(\dom R^{\widetilde{*}}\bigr)\) is closed in \(\widehat{\mathfrak H}\). The following statements hold.
\begin{enumerate}
\renewcommand*\theenumi{\alph{enumi}}
\renewcommand*\labelenumi{\rm(\theenumi)}
\item  \label{th-comadj4-i1}
If \(R\) is closed, that is, \(R = \widebar{R}\), then \(\CH(R)\) is closed.
\item \label{th-comadj4-i2}
If \(R\) is maximal dissipative  in \( \widetilde{\mathfrak{H}}\),  then \(\CH(R)\) is a maximal dissipative linear relation in \(\mathfrak H\).
\item \label{th-comadj4-i3}
If \(R\) is maximal symmetric in \(\widetilde{\mathfrak{H}}\),  then \(\CH(R)\) is a maximal symmetric linear relation in \(\mathfrak H\).
\item \label{th-comadj4-i4}
If \(R\) is self-adjoint in \(\widetilde{\mathfrak{H}}\), then \(\CH(R)\) is a self-adjoint relation in \(\mathfrak{H}\).
\end{enumerate}
\end{theorem}

\begin{proof}
Let \(R\) be a linear relation in \(\widetilde{\mathfrak{H}}\) and assume that  \(\widetilde{P}_{\widehat{\mathfrak H}}\bigl(\dom R^{\widetilde{*}}\bigr)\) is closed in \(\widehat{\mathfrak H}\).

Item (\ref{th-comadj4-i1}) follows from Corollary~\ref{co-thR*}.

(\ref{th-comadj4-i2})
Assume \(R\) is a maximal dissipative linear relation in \(\bigl(\widetilde{\mathfrak{H}}, \kip_{\widetilde{\mathfrak{H}}} \bigr)\). Equivalently, by Theorem~\ref{th-maxnnR}(\ref{th-maxnnR-i1}), assume \(R\) is closed and both \(R\) and \(-R^{\widetilde{*}}\) are dissipative. We prove (\ref{th-comadj4-i2}) by showing that these three properties also hold for \(\CH(R)\) in place of \(R\). By (\ref{th-comadj4-i1}), \(\CH(R)\) is closed and by
Proposition~\ref{pr-inher}(\ref{pr-inher-i2}), \(\CH(R)\) and \(\CH(-R^{\widetilde{*}})\) are dissipative. Corollary~\ref{co-thR*} implies
\[
-\CH(R)^*
=
-\overline{\CH(R^{\widetilde{*}})}
=
\overline{\CH(-R^{\widetilde{*}})}. 
\]
Since the closure of a dissipative relation is dissipative, \(-\CH(R)^*\) is dissipative. Hence, by Theorem~\ref{th-maxnnR}(\ref{th-maxnnR-i1}),
\(\CH(R)\) is maximal dissipative.

(\ref{th-comadj4-i3})
Assume that \(R\) is a maximal symmetric linear relation in \(\bigl( \widetilde{\mathfrak{H}}, \kip_{\widetilde{\mathfrak{H}}} \bigr)\). By Theorem~\ref{th-maxnnR}(\ref{th-maxnnR-i2}), \(R\) is symmetric and \(R\) or \(-R\) is maximal dissipative. By Proposition~\ref{pr-inher}(\ref{pr-inher-i2}), \(\CH(R)\) is symmetric. Moreover, applying (\ref{th-comadj4-i2})  to \(R\) or to \(-R\) shows that \(\CH(R)\) or \(-\CH(R)\) is maximal dissipative. Again, by Theorem~\ref{th-maxnnR}(\ref{th-maxnnR-i2}), \(\CH(R)\) is maximal symmetric.

(\ref{th-comadj4-i4})
Assume \(R = R^{\widetilde{*}}\). Then \(R\) is closed, and by the main assumption \(\widetilde{P}_{\widehat{\mathfrak H}}(\dom R)=\widetilde{P}_{\widehat{\mathfrak H}}\bigl(\dom R^{\widetilde{*}}\bigr)\) is closed. By Theorem~\ref{th-comadjn}, \(\CH(R) = \CH(R^{\widetilde{*}}) = \CH(R)^*\). Hence \(\CH(R)\) is self-adjoint, proving (\ref{th-comadj4-i4}).
\end{proof}

In the next theorem we give reflection results, that is, converse-type statements to the corresponding preservation results in Theorem~\ref{th-comadj4}. These reverse implications require not only the closedness of \(\widetilde P_{\widehat{\mathfrak H}}(\dom R)\) in \(\widehat{\mathfrak H}\), but also the equality of relevant  multivalued parts. In the proof we use the following lemma.

\begin{lemma}\label{le-R=S}
Let \(\widetilde{\mathfrak H} = \mathfrak{H} [\dot{+}] \widehat{\mathfrak H}\) be a Krein space and let \(R\) and \(S\) be linear relations in \(\widetilde{\mathfrak H}\) such that  \(R \subseteq S\). If \(S(0) = R(0)\), \(\widetilde P_{\widehat{\mathfrak H}}(\dom S) = \widetilde P_{\widehat{\mathfrak H}}(\dom R)\) and  \(\CH(S)=\CH(R)\),
then \(S=R\).
\end{lemma}

\begin{proof}
Let \(\{\widetilde x,\widetilde y\} \in S\). By the equality of the projected
domains, there exists \(\{\widetilde f,\widetilde g\}\in R\) such that
\(\widetilde P_{\widehat{\mathfrak H}}\widetilde f = \widetilde P_{\widehat{\mathfrak H}}\widetilde x\). Since \(R \subseteq S\), we have \(\{\widetilde x - \widetilde f,\widetilde y - \widetilde g\}\in S\), and \(\widetilde x - \widetilde f \in \mathfrak H\). Set \(\{\widetilde x - \widetilde f,\widetilde y - \widetilde g\} = \{h,\widetilde{v}\} \in S|_{\mathfrak H}\). Since \(\CH(S)=\CH(R)\), there exists \(\widetilde{k} \in \widetilde{\mathfrak{H}}\) such that \(\{h,\widetilde{k}\} \in R|_{\mathfrak H}\) and \(\widetilde P_{\mathfrak H}\widetilde{k} = \widetilde P_{\mathfrak H}\widetilde{v}\). Now, \(R|_{\mathfrak H} \subseteq S|_{\mathfrak H}\) yields \(\{0,\widetilde{v} - \widetilde{k}\} \in S\), and \(S(0) = R(0)\) implies \(\{0,\widetilde{v} - \widetilde{k}\} \in R\). Thus, \(\{\widetilde f,\widetilde g\}, \{h,\widetilde{k}\}, \{0,\widetilde{v} - \widetilde{k}\} \in R\) and
\[
\{\widetilde x,\widetilde y\} = \{\widetilde f,\widetilde g\} +
\{h,\widetilde{k}\}  + \{0,\widetilde{v} - \widetilde{k}\}  \in R,
\]
proving \(S\subseteq R\). Since by assumption \(R\subseteq S\), the lemma is proved.
\end{proof}

\begin{theorem}\label{th-conv}
Let \(\widetilde{\mathfrak H}=\mathfrak H[\dotplus]\widehat{\mathfrak H}\)
be a Krein space and let \(R\) be a linear relation in \(\widetilde{\mathfrak H}\). Assume \(\widetilde{P}_{\widehat{\mathfrak H}}(\dom R)\) is closed in \(\widehat{\mathfrak H}\). The following statements hold.

\begin{enumerate}
\renewcommand*\theenumi{\alph{enumi}}
\renewcommand*\labelenumi{\rm(\theenumi)}

\item \label{th-conv-i1}
If \(R(0)=\widebar R(0)\) and \(\CH(R)=\CH(\widebar R)\),
then \(R\) is closed, that is, \(R= \widebar R\).

\item \label{th-conv-i2}
If \(R\) is dissipative, \(R(0)=R^{\widetilde{*}}\mkern-2mu(0)\), and \(\CH(R)\) is maximal dissipative in \(\mathfrak H\), then \(R\) is maximal
dissipative in \(\widetilde{\mathfrak H}\).

\item \label{th-conv-i3}
If \(R\) is symmetric, \(R(0)=R^{\widetilde{*}}\mkern-2mu(0)\), and \(\CH(R)\) is maximal symmetric in \(\mathfrak H\), then \(R\) is maximal symmetric in \(\widetilde{\mathfrak H}\).

\item \label{th-conv-i4}
If \(R\) is symmetric, \(R(0)=R^{\widetilde{*}}\mkern-2mu(0)\), and \(\CH(R)\) is self-adjoint in \(\mathfrak H\), then \(R\) is self-adjoint in \(\widetilde{\mathfrak H}\).

\end{enumerate}
\end{theorem}

\begin{proof}
(\ref{th-conv-i1})
To apply Lemma~\ref{le-R=S} with \(S = \widebar R\) we need to prove
\begin{equation}\label{eq-pdR-pdbR}
\widetilde P_{\widehat{\mathfrak H}}(\dom R)
=
\widetilde P_{\widehat{\mathfrak H}}\bigl(\dom\widebar R\bigr).
\end{equation}
Since the inclusion \(\subseteq\) is trivial, we prove \(\supseteq\) next. Let \(\widehat y \in \widetilde P_{\widehat{\mathfrak H}}(\dom\widebar R)\) be arbitrary. Then there exists \(\{\widetilde f,\widetilde g\}\in\widebar R\) such that \(\widehat y = \widetilde P_{\widehat{\mathfrak H}}\widetilde f\). As \(\{\widetilde f,\widetilde g\}\in\widebar R\), there exists a sequence \(\{\widetilde f_n,\widetilde g_n\}\in R\) such that
\(\{\widetilde f_n,\widetilde g_n\}\to \{\widetilde f,\widetilde g\}\) as \(n\to\infty\).
Hence
\(\widetilde P_{\widehat{\mathfrak H}}\widetilde f_n
\to
\widetilde P_{\widehat{\mathfrak H}}\widetilde f = \widehat y\) as \(n\to\infty\).
Since \(\widetilde P_{\widehat{\mathfrak H}}\widetilde f_n \in \widetilde P_{\widehat{\mathfrak H}}(\dom R)\) and \(\widetilde P_{\widehat{\mathfrak H}}(\dom R)\) is closed, we have \(\widehat y \in \widetilde P_{\widehat{\mathfrak H}}(\dom R)\) proving \(\supseteq\) in \eqref{eq-pdR-pdbR}. Now Lemma~\ref{le-R=S} yields \(R = \widebar R\), that is, \(R\) is closed.

(\ref{th-conv-i2}) Assume that \(R\) is dissipative, \(R(0)=R^{\widetilde{*}}\mkern-2mu(0)\), \(\widetilde P_{\widehat{\mathfrak H}}(\dom R)\) is closed, and \(\CH(R)\) is maximal dissipative. Let \(S\) be any dissipative linear relation such that \(R\subseteq S\). Lemma~\ref{lem-mdd} and reversal of inclusion under adjoints give
\[
R(0)\subseteq S(0)\subseteq S^{\widetilde{*}}\mkern-2.5mu(0)
\subseteq R^{\widetilde{*}}\mkern-2mu(0)=R(0).
\]
Thus \(S(0)=R(0)\) and \(S^{\widetilde{*}}\mkern-2mu(0)=R^{\widetilde{*}}\mkern-2mu(0)\). Equality \eqref{eq-R*(0)} yields \(\overline{\dom S}=\overline{\dom R}\). Continuity of \(\widetilde P_{\widehat{\mathfrak H}}\) and closedness of \(\widetilde P_{\widehat{\mathfrak H}}(\dom R)\) imply
\[
\widetilde P_{\widehat{\mathfrak H}}\bigl(\dom S\bigr)\subseteq\widetilde P_{\widehat{\mathfrak H}}\bigl(\overline{\dom R}\bigr)
\subseteq\overline{\widetilde P_{\widehat{\mathfrak H}}(\dom R)}=\widetilde P_{\widehat{\mathfrak H}}(\dom R).
\]
The reverse inclusion follows from \(R\subseteq S\), hence \(\widetilde P_{\widehat{\mathfrak H}}\bigl(\dom S\bigr) = \widetilde P_{\widehat{\mathfrak H}}(\dom R)\). Finally, \(\CH(S)\) is dissipative and contains \(\CH(R)\), so maximality gives \(\CH(S)=\CH(R)\). Lemma~\ref{le-R=S} now gives \(S=R\). Hence \(R\) is maximal dissipative.

(\ref{th-conv-i3})
Assume that \(R\) is symmetric, \(R(0)=R^{\widetilde{*}}\mkern-2mu(0)\), and \(\CH(R)\) is maximal symmetric. By Theorem~\ref{th-maxnnR}{(\ref{th-maxnnR-i2})}, 
\(\CH(R)\) or \(-\CH(R)\) is maximal dissipative.

We have \(\CH(-R)=-\CH(R)\), \(\widetilde P_{\widehat{\mathfrak H}}\bigl(\dom(-R)\bigr)
= \widetilde P_{\widehat{\mathfrak H}}(\dom R)\), and \((-R)^{\widetilde{*}}=-R^{\widetilde{*}}\). Also, since \(R\) is symmetric, both \(R\) and \(-R\) are dissipative, and \((-R)(0)=(-R)^{\widetilde{*}}\mkern-2mu(0)\) follows from \(R(0)=R^{\widetilde{*}}\mkern-2mu(0)\). Hence part~(\ref{th-conv-i2}), applied to \(R\) or to \(-R\), gives that \(R\) or \(-R\) is maximal dissipative. Since \(R\) is symmetric, Theorem~\ref{th-maxnnR}{\rm(\ref{th-maxnnR-i2})} implies that \(R\) is maximal symmetric.

(\ref{th-conv-i4})
Assume that \(R\) is symmetric, \(R(0)=R^{\widetilde{*}}\mkern-2mu(0)\), and \(\CH(R)\) is self-adjoint. By Theorem~\ref{th-maxnnR}{\rm(\ref{th-maxnnR-i3})}, both \(\CH(R)\) and \(-\CH(R)=\CH(-R)\) are maximal dissipative. As observed in the proof of {\rm(\ref{th-conv-i3})}, the hypotheses of part~{\rm(\ref{th-conv-i2})} are satisfied for both \(R\) and \(-R\). Therefore \(R\) and \(-R\) are maximal dissipative. By Theorem~\ref{th-maxnnR}{\rm(\ref{th-maxnnR-i3})}, \(R\) is self-adjoint.
\end{proof}

The preservation Theorem~\ref{th-comadj4} could be stated under the hypothesis \(\CH(\widebar R)=\CH(R^{\widetilde{*}})^*\) which, by Corollary~\ref{co-thR*}, is weaker than the hypothesis that \(\widetilde P_{\widehat{\mathfrak H}}\bigl(\dom R^{\widetilde{*}}\bigr)\) is closed. However, the proof of the reflection Theorem~\ref{th-conv} essentially uses the closedness of \(\widetilde{P}_{\widehat{\mathfrak H}}(\dom R)\). We have stated Theorem~\ref{th-comadj4} with the stronger projected-domain hypothesis to keep the preservation and reflection results in parallel form. The following corollary makes this parallelism explicit by combining Theorems~\ref{th-comadj4} and~\ref{th-conv} with Lemma~\ref{lem-mdd}.

\begin{corollary}\label{th-4-iff}
Let \(\widetilde{\mathfrak H}=\mathfrak H[\dotplus]\widehat{\mathfrak H}\) be a Krein space, and let \(R\) be a linear relation in \(\widetilde{\mathfrak H}\). Assume that both subspaces \(\widetilde P_{\widehat{\mathfrak H}}(\dom R)\) and \(\widetilde P_{\widehat{\mathfrak H}}\bigl(\dom R^{\widetilde{*}}\bigr)\) are closed in \(\widehat{\mathfrak H}\). The following statements hold.

\begin{enumerate}
\renewcommand*\theenumi{\alph{enumi}}
\renewcommand*\labelenumi{\rm(\theenumi)}

\item \label{th-4-iff-i1}
\(R\) is closed, that is, \(R = \widebar{R}\),  if and only if \(R(0)=\widebar R(0)\) and \(\CH(R)=\CH(\widebar R)\).

\item \label{th-4-iff-i2}
\(R\) is maximal dissipative if and only if \(R\) is dissipative, \(R(0)=R^{\widetilde{*}}\mkern-2mu(0)\), and \(\CH(R)\) is maximal dissipative in \(\mathfrak H\).

\item \label{th-4-iff-i3}
\(R\) is maximal symmetric if and only if \(R\) is symmetric, \(R(0)=R^{\widetilde{*}}\mkern-2mu(0)\) and \(\CH(R)\) is maximal symmetric in \(\mathfrak H\).

\item \label{th-4-iff-i4}
\(R\) is self-adjoint if and only if \(R\) is symmetric, \(R(0)=R^{\widetilde{*}}\mkern-2mu(0)\) and \(\CH(R)\) is self-adjoint in \(\mathfrak H\).

\end{enumerate}
\end{corollary}

\begin{proof}
The ``only if'' part of (\ref{th-4-iff-i1}) is trivial, while the ``if'' part is (\ref{th-conv-i1}) in Theorem~\ref{th-conv}.

(\ref{th-4-iff-i2}) Assume that \(R\) is maximal dissipative. Then \(R\) is clearly dissipative and it is closed by Theorem~\ref{th-maxnnR}(\ref{th-maxnnR-i1}). By Lemma~\ref{lem-mdd} we have \(R(0)=R^{\widetilde{*}}\mkern-2mu(0)\) and \(\CH(R)\) is maximal dissipative in \(\mathfrak H\) by Theorem~\ref{th-comadj4}(\ref{th-comadj4-i2}). The converse is (\ref{th-conv-i2}) in Theorem~\ref{th-conv}.

(\ref{th-4-iff-i3}) Assume that \(R\) is maximal symmetric. Then \(R\) is clearly symmetric and it is closed since the closure of a neutral subspace of the Lagrange-Krein space \(\bigl(\widetilde{\mathfrak{H}}^2, \lfkip\bigr)\) is neutral in this space. By Lemma~\ref{lem-mdd} we have \(R(0)=R^{\widetilde{*}}\mkern-2mu(0)\) and \(\CH(R)\) is maximal symmetric in \(\mathfrak H\) by Theorem~\ref{th-comadj4}(\ref{th-comadj4-i3}). The converse is (\ref{th-conv-i3}) in Theorem~\ref{th-conv}.

(\ref{th-4-iff-i4}) Assume that \(R\) is self-adjoint. Then, clearly, \(R\) is  symmetric and closed. The equality \(R(0)=R^{\widetilde{*}}\mkern-2mu(0)\) is trivial and \(\CH(R)\) is self-adjoint in \(\mathfrak H\) by Theorem~\ref{th-comadj4}(\ref{th-comadj4-i4}). The converse is (\ref{th-conv-i4}) in Theorem~\ref{th-conv}.
\end{proof}

\section{Compression of the Friedrichs extension} \label{sec-Friedrichs}

A linear relation \(S\) in a Krein space \(\bigl(\mathfrak{K},\kip_{\mathfrak{K}}\bigr)\) is called \emph{nonnegative} if \([y,x]_{\mathfrak{K}} \geq 0\) for all \(\{x,y\}\in S\).

For a nonnegative linear relation \(S\)  we have
\([y,x]_{\mathfrak{K}} \in \mathbb{R}\) for all \(\{x,y\}\in S\).
The characterization of symmetry mentioned after \eqref{eq-dissym}
shows that every nonnegative linear relation is
symmetric. We therefore omit the adjective ``symmetric'' when referring to nonnegative linear relations.  

Nonnegativity of a linear relation \(S\) should be distinguished from nonnegativity of \(S\) as a subspace of the Lagrange-Krein space. As the table in Subsection~\ref{subs-ddip} shows, the latter corresponds to dissipativity of the linear relation. Moreover, neutrality of the subspace implies symmetry of the linear relation, but does not imply nonnegativity of the relation.

\begin{proposition} \label{pr-inhnn}
Let \(\widetilde{\mathfrak H}=\mathfrak H [\dotplus] \widehat{\mathfrak H}\) be a Krein space and let \(S\) be a linear relation in \(\widetilde{\mathfrak H}\). If \(S\) is nonnegative in \(\bigl(\widetilde{\mathfrak H},\kip_{\widetilde{\mathfrak H}}\bigr)\), then \(\CH(S)\) is nonnegative in \(\bigl(\mathfrak H,\kip_{\mathfrak H}\bigr)\).
\end{proposition}

\begin{proof}
Assume that \(S\) is nonnegative. Let \(\{f,h\}\in\CH(S)\). Then there exists \(\{f,\widetilde g\}\in S|_{\mathfrak H}\) such that \(h=\widetilde P_{\mathfrak H}\widetilde g\). Since \(f\in\mathfrak H\), we have
\[
[h,f]_{\mathfrak H}
=
[\widetilde P_{\mathfrak H}\widetilde g,f]_{\mathfrak H}
=
[\widetilde g,f]_{\widetilde{\mathfrak H}}.
\]
The last expression is nonnegative because \(S\) is nonnegative. Hence
\(\CH(S)\) is nonnegative in \(\mathfrak H\). This proves the proposition. 
\end{proof}

The following definition extends the Hilbert space definition
in \cite[p.~207]{CdS1978} to Krein spaces. 

\begin{definition} \label{de-Fe} 
Let \(S\) be a nonnegative linear relation in the Krein space
\(\bigl(\mathfrak{K}, \kip_{\mathfrak K}\bigr)\). 
The \emph{Friedrichs extension} of \(S\), denoted by \(S_{\mathrm F}\), consists of all pairs \(\{u,v\}\in S^{[*]}\) for which there exists a sequence \(\bigl(\{f_n,g_n\}\bigr)_{n\in\mathbb N}\) in \(S\) such that \(f_n\to u\) in \(\mathfrak{K}\) and \(\bigl[g_n-g_m,f_n-f_m\bigr]_{\mathfrak{K}}\to 0\) as \(n,m\to\infty\). Here \(S^{[*]}\) denotes the adjoint of \(S\) in the Krein space \(\bigl(\mathfrak{K}, \kip_{\mathfrak K}\bigr)\). 
\end{definition}

\begin{proposition}
Let \(\bigl(\mathfrak{K}, \kip_{\mathfrak K}\bigr)\) be a Krein space, let \(J\) be a fundamental symmetry on \(\mathfrak{K}\), and let \(\ahip_{\mathfrak{K}} = [J \mkern1mu\cdot\mkern1mu,\mkern-2mu \cdot\mkern1.5mu ]_{\mathfrak{K}}\) 
be the corresponding Hilbert space inner product. Let \(S\) be a linear relation in \(\mathfrak{K}\). Then \(S\) is a nonnegative linear relation in the Krein space \(\mathfrak{K}\) if and only if \(JS\) is a nonnegative linear relation in the Hilbert space \(\bigl(\mathfrak{K},\ahip_{\mathfrak{K}} \bigr)\).  If these equivalent conditions hold, then \(S_{\mathrm F}=J((JS)_{\mathrm F})\). Hence, \(S_{\mathrm F}\) is a nonnegative self-adjoint extension of \(S\) in \(\bigl(\mathfrak{K}, \kip_{\mathfrak K}\bigr)\). 
Here \((JS)_{\mathrm F}\) is the Friedrichs extension of \(JS\) in the Hilbert space \(\bigl(\mathfrak{K},\ahip_{\mathfrak{K}} \bigr)\).
\end{proposition}

\begin{proof}
We denote the adjoint in the Hilbert space \(\bigl(\mathfrak{K},\ahip_{\mathfrak{K}}\bigr)\) by \(^{\langle *\rangle}\). Let \(\{u,v\} \in \mathfrak{K}^2\) be arbitrary. The following equivalences hold 
\begin{align*}
\{u,Jv\} \in J S^{[*]} 
& \ \ \Leftrightarrow \ \ 
\forall \mkern2mu \{x,y\} \in S \quad [y,u]_{\mathfrak{K}} = [x,v]_{\mathfrak{K}} \\ 
& \ \ \Leftrightarrow \ \ 
\forall \mkern2mu \{x,Jy\} \in JS \quad \langle Jy,u\rangle_{\mathfrak{K}} = \langle x, J v\rangle_{\mathfrak{K}} \\   
& \ \ \Leftrightarrow \ \ 
\{u,Jv\} \in (J S)^{\langle *\rangle}. 
\end{align*}
Hence \(J S^{[*]} = (JS)^{\langle *\rangle}\). 

For all \(\{f,g\} \in S\) we have \([g,f]_{\mathfrak{K}} = \langle Jg, f \rangle_{\mathfrak{K}}\). Thus \(S\) is nonnegative in \(\bigl(\mathfrak{K}, \kip_{\mathfrak K}\bigr)\) if and only if \(JS\) is nonnegative in \(\bigl(\mathfrak{K}, \ahip_{\mathfrak K}\bigr)\).

Assume that \(S\) is nonnegative in \(\bigl(\mathfrak{K}, \kip_{\mathfrak K}\bigr)\), or equivalently, \(JS\) is nonnegative in \(\bigl(\mathfrak{K}, \ahip_{\mathfrak K}\bigr)\). Let \((JS)_{\mathrm F}\) be the Friedrichs extension of \(JS\) in the Hilbert space \(\bigl(\mathfrak{K}, \ahip_{\mathfrak K}\bigr)\) as defined in \cite[Section~3]{CdS1978}. Then \((JS)_{\mathrm F}\) is a nonnegative self-adjoint extension of \(JS\) in \(\bigl(\mathfrak{K}, \ahip_{\mathfrak K}\bigr)\). It consists of all pairs \(\{x,y\}\in (JS)^{\langle *\rangle}\) for which there exists a sequence \(\bigl(\{f_n,Jg_n\}\bigr)_{n\in\mathbb N}\) in \(JS\) such that \(f_n\to x\) in \(\mathfrak{K}\) and \(\bigl\langle Jg_n - Jg_m, f_n-f_m \bigr\rangle_{\mathfrak{K}} \to 0\) as \(n,m\to\infty\). 
Since \((JS)^{\langle *\rangle}=JS^{[*]}\) and
\[
\bigl\langle Jg_n-Jg_m,f_n-f_m \bigr\rangle_{\mathfrak K}
=
\bigl[g_n-g_m,f_n-f_m \bigr]_{\mathfrak K},
\]
the set described in the preceding sentence is precisely \(JS_{\mathrm F}\). Hence \((JS)_{\mathrm F}=JS_{\mathrm F}\). Applying \(J\) to both sides gives \(S_{\mathrm F} = J\bigl((JS)_{\mathrm F}\bigr)\). Consequently, \(S_{\mathrm F}\) is a nonnegative self-adjoint extension of \(S\) in \(\bigl(\mathfrak{K}, \kip_{\mathfrak K}\bigr)\).
\end{proof}

\begin{theorem} \label{th-FriCom}
Let \(\widetilde{\mathfrak H}=\mathfrak H [\dotplus] \widehat{\mathfrak H}\) be a Krein space and let \(S\) be a nonnegative linear relation in \(\widetilde{\mathfrak H}\). If \(\CH(S^{\widetilde{*}}) = \CH(S)^*\), then \(\CH(S_{\mathrm F})= (\CH(S))_{\mathrm F}\). 
\end{theorem}

\begin{proof}
By Proposition~\ref{pr-inhnn}, the compression \(\CH(S)\) is a nonnegative
linear relation in \(\mathfrak H\).

Assume that \(\CH(S^{\widetilde{*}})=\CH(S)^*\). We first prove that \( (\CH(S) )_{\mathrm F}\subseteq \CH(S_{\mathrm F})\). Let \(\{f,g\}\in (\CH(S))_{\mathrm F}\). Then \(\{f,g\}\in\CH(S)^*=\CH(S^{\widetilde{*}})\). Hence there exists \(\widetilde g\in\widetilde{\mathfrak H}\) such that \(\{f,\widetilde g\}\in S^{\widetilde{*}}\) and \(\widetilde P_{\mathfrak H}\widetilde g=g\). By the definition of the Friedrichs extension, there exists a sequence \(\{f_n,g_n\}\in \CH(S)\) such that 
\[
f_n\to f \quad \text{and} \quad \bigl[g_n-g_m,f_n-f_m \bigr]_{\mathfrak H}\to 0
\quad \text{as} \quad n,m\to\infty .
\]
For every \(n\), choose \(\widetilde g_n\in\widetilde{\mathfrak H}\) such that
\(\{f_n,\widetilde g_n\}\in S\) and  \(\widetilde P_{\mathfrak H}\widetilde g_n=g_n\). 
Since \(f_n-f_m\in\mathfrak H\), we have
\[
\bigl[
\widetilde g_n-\widetilde g_m,f_n-f_m \bigr]_{\widetilde{\mathfrak H}}
=
\bigl[
g_n-g_m,f_n-f_m\bigr]_{\mathfrak H} \to 0 \quad \text{as} \quad n,m\to\infty .
\]
Therefore, by the definition of \(S_{\mathrm F}\) we have \(\{f,\widetilde g\}\in S_{\mathrm F}\). It follows that
\[
\{f,g\}
=
\{f,\widetilde P_{\mathfrak H}\widetilde g\}
\in \CH(S_{\mathrm F}).
\]
Thus
\begin{equation} \label{eq-inclCH}
 (\CH(S) )_{\mathrm F} \subseteq \CH(S_{\mathrm F}).
\end{equation} 

Since \(S_{\mathrm F}\) is self-adjoint, by Proposition~\ref{pr-inher}(\ref{pr-inher-i2}), \(\CH(S_{\mathrm F})\) is symmetric. Since a self-adjoint relation has no proper symmetric extensions, the inclusion in \eqref{eq-inclCH} implies the desired equality. 
\end{proof}

The preceding proof follows the lifting argument in \cite[proof of Theorem~4.5]{ADW2013}. There, finite codimension is used to establish self-adjointness of \(\CH(S_{\mathrm F})\) before proving the inclusion \eqref{eq-inclCH}. Here symmetry of \(\CH(S_{\mathrm F})\) suffices, since \((\CH(S))_{\mathrm F}\) is self-adjoint. Thus the argument also applies to infinite-dimensional exit spaces whenever the adjoint identity holds. 

\begin{remark}
There is an analogous question for the other extremal nonnegative
self-adjoint extension. In Hilbert space terminology this extension is often
called the Krein--von Neumann extension; in Krein's terminology it is the soft
extension. Its behavior under compression is more delicate than that
of the Friedrichs extension. For finite-codimensional compressions in Hilbert
spaces, this question is studied in \cite{ADW2013}. As \cite[Example~4.6]{ADW2013} shows, the identity \(\CH(S^{\widetilde{*}})=\CH(S)^*\) is not sufficient to obtain the corresponding compression formula for this extension. The positive result in \cite[Theorem~4.7]{ADW2013} requires an additional range condition. We do not pursue this direction here. 
\end{remark}

\section{Proof of Theorem~\ref{equivalences}} \label{SecEqi}

\begin{lemma} \label{le-nidom}
In an infinite-dimensional Krein space \(\bigl(\mathfrak{K},\kip_{\mathfrak{K}}\bigr)\) for every vector \(v \in \mathfrak{K}\setminus\{0\}\) there exists a densely defined self-adjoint operator \(B\) in \(\mathfrak{K}\) such that \(v \not\in \dom B\).
\end{lemma}

\begin{proof}
Let \(J\) be a fundamental symmetry on \(\bigl(\mathfrak{K},\kip_{\mathfrak{K}}\bigr)\) and let \(\ahip_{\mathfrak{K}}\) be the corresponding Hilbert space inner product. Since \(\mathfrak K\) is infinite dimensional, there exists a densely defined self-adjoint operator \(A\) in the Hilbert space \(\bigl(\mathfrak{K},\ahip_{\mathfrak{K}}\bigr)\) with \(\dom A\neq\mathfrak K\). For example, one may take a diagonal unbounded self-adjoint operator on a separable infinite-dimensional closed subspace and extend it by zero on the orthogonal complement.

If \(v \notin \dom A\), then \(B=JA\) is the desired operator since \(B\) is self-adjoint in the Krein space \(\bigl(\mathfrak K,\kip_{\mathfrak K}\bigr)\) as \(JB=A\) is self-adjoint in the Hilbert space \(\bigl(\mathfrak K,\ahip_{\mathfrak K}\bigr)\).

If \(v \in \dom A\), then choose \(u \not\in \dom A\). Since replacing the prescribed vector by its normalization does not change whether it belongs to an operator domain, we can choose \(u\) and \(v\) to be unit vectors and \(\langle u,v  \rangle_{\mathfrak{K}} = 0\). Indeed, if \(\langle u,v  \rangle_{\mkern-1mu\mathfrak{K}} \neq 0\), then \(u' = \langle u,v\rangle_{\mkern-2mu\mathfrak{K}}\mkern1mu v - u\) is orthogonal to \(v\) and \(u' \not\in \dom A\), because \(v\in\dom A\) while \(u\notin\dom A\). Replacing \(u\) by the normalization of \(u'\), we may therefore assume \(\langle u,v\rangle_{\mathfrak K}=0\).

Define \(U:\mathfrak K\to\mathfrak K\) by
\(
U(\alpha u+\beta v)=-\alpha v+\beta u
\)
for all \( \alpha,\beta\in\mathbb C\), and let \(Uw=w\) for all \(w\perp\lspan\{u,v\}\) in \(\bigl(\mathfrak K,\ahip_{\mathfrak K}\bigr)\). Then \(U\) is unitary and \(Uv=u\).
Set \(B=JU^{-1}AU\). Then \(JB=U^{-1}AU\) is self-adjoint in the Hilbert space \(\bigl(\mathfrak K,\ahip_{\mathfrak K}\bigr)\), and hence \(B\) is self-adjoint in the Krein space \(\bigl(\mathfrak K,\kip_{\mathfrak K}\bigr)\). Moreover, \(\dom B=U^{-1}(\dom A)\), and since \(Uv=u\notin\dom A\), we have \(v\notin\dom B\).
\end{proof}

\begin{lemma} \label{le-comcl}
Let \(\widetilde{\mathfrak H}=\mathfrak H [\dotplus] \widehat{\mathfrak H}\) be a Krein space and assume that \(\widehat{\mathfrak H}\) is finite dimensional. Then \(\CH(\widebar R)=\overline{\CH(R)}\) for every linear relation \(R\) in \(\widetilde{\mathfrak H}\).
\end{lemma}
\begin{proof}
Since \(\widehat{\mathfrak H}\) is finite dimensional, its subspace
\(\widetilde{P}_{\widehat{\mathfrak H}}(\operatorname{dom} R^{\widetilde{*}})\)
is closed. By
Corollary~\ref{co-thR*}, the relation \(\CH(\widebar{R})\) is closed.
Consequently, the inclusion \(\CH(R)\subseteq\CH(\widebar{R})\) implies that
\[
\overline{\CH(R)}\subseteq\CH(\widebar{R}).
\]

To prove the reverse inclusion, let
\(\{f,\widetilde{g}\}\in\widebar{R}\) with \(f\in\mathfrak H\).
Choose a sequence
\((\{\widetilde{x}_n,\widetilde{y}_n\})_{n\in\mathbb N}\)
in \(R\) such that
\[
\{\widetilde{x}_n,\widetilde{y}_n\}
\to \{f,\widetilde{g}\}
\quad \text{as} \quad n\to\infty.
\]
Here and below, convergence of pairs is with respect to the Hilbert
space topology of
\(\bigl(\widetilde{\mathfrak H}^{2}, \dahip\bigr)\).
We will modify these pairs so that their first components belong
to \(\mathfrak H\), without changing their limit.

Define the linear mapping
\[
Q:R\to
\widetilde{P}_{\widehat{\mathfrak H}}(\operatorname{dom} R)
\]
by
\[
Q\{\widetilde{x},\widetilde{y}\}
=
\widetilde{P}_{\widehat{\mathfrak H}}\widetilde{x}, \quad \{\widetilde{x},\widetilde{y}\}\in R.
\]
Both its domain and range carry the topologies inherited from
the corresponding Hilbert spaces. By its definition, \(Q\) is continuous and surjective, and \(\operatorname{nul} Q=R|_{\mathfrak H}\). Moreover, continuity of
\(\widetilde{P}_{\widehat{\mathfrak H}}\) and the fact that
\(f\in\mathfrak H\) give
\[
Q\{\widetilde{x}_n,\widetilde{y}_n\}
=
\widetilde{P}_{\widehat{\mathfrak H}}\widetilde{x}_n
\to
\widetilde{P}_{\widehat{\mathfrak H}}f
=
0
\quad \text{as} \quad n\to\infty.
\]

We next construct a finite-dimensional complement \(\mathfrak{M}\) of
\(R|_{\mathfrak H}\) in \(R\). Let
\[
j=
\dim\bigl(
\widetilde{P}_{\widehat{\mathfrak H}}(\operatorname{dom} R)
\bigr).
\]
If \(j=0\), set \(\mathfrak{M} = \bigl\{\{0,0\} \bigr\}\). If \(j>0\), let
\(\{\widehat{z}_1,\ldots,\widehat{z}_j\}\) be a basis for
\(\widetilde{P}_{\widehat{\mathfrak H}}(\operatorname{dom} R)\).
Since \(Q\) is surjective, we can choose
\(\{\widetilde{f}_i,\widetilde{g}_i\}\in R\) such that
\[
Q\{\widetilde{f}_i,\widetilde{g}_i\}
=
\widehat{z}_i
\quad \text{for all} \quad i\in\{1,\ldots,j\}.
\]
Set
\[
\mathfrak{M} =
\operatorname{span}
\bigl\{
\{\widetilde{f}_1,\widetilde{g}_1\},
\ldots,
\{\widetilde{f}_j,\widetilde{g}_j\}
\bigr\}.
\]
Applying \(Q\) to a linear combination of these pairs shows that
they are linearly independent, because their images form a basis.
Thus \(\mathfrak{M}\) is \(j\)-dimensional, and
\[
Q|_\mathfrak{M} : \mathfrak{M} \to
\widetilde{P}_{\widehat{\mathfrak H}}(\operatorname{dom} R)
\]
is a linear bijection. This conclusion also holds when \(j=0\).

For every pair in \(R\), there is therefore a unique pair in \(\mathfrak{M}\)
with the same image under \(Q\). Their difference belongs to
\(\operatorname{nul} Q=R|_{\mathfrak H}\). Since \(Q|_\mathfrak{M}\) is
injective, \(\mathfrak{M}\cap R|_{\mathfrak H}=\bigl\{\{0,0\}\bigr\}\), and hence
\[
R=\mathfrak{M}\dotplus R|_{\mathfrak H}.
\]
Moreover, since both the domain and the codomain of \(Q|_{\mathfrak{M}}\)
are finite dimensional, its inverse
\[
(Q|_\mathfrak{M})^{-1}:
\widetilde{P}_{\widehat{\mathfrak H}}(\operatorname{dom} R)
\to \mathfrak{M}
\]
is continuous.

Using this direct sum decomposition, write, for every
\(n\in\mathbb N\),
\[
\{\widetilde{x}_n,\widetilde{y}_n\}
=
\{\widetilde{h}_n,\widetilde{k}_n\}
+
\{u_n,\widetilde{v}_n\},
\]
where
\(\{\widetilde{h}_n,\widetilde{k}_n\}\in \mathfrak{M}\) and
\(\{u_n,\widetilde{v}_n\}\in R|_{\mathfrak H}\).
Since \(Q\) vanishes on \(R|_{\mathfrak H}\), we have
\[
Q|_{\mathfrak{M}}\{\widetilde{h}_n,\widetilde{k}_n\}
=
Q\{\widetilde{x}_n,\widetilde{y}_n\}
\to 0
\quad \text{as} \quad n\to\infty.
\]
Continuity of \((Q|_{\mathfrak{M}})^{-1}\) now gives
\[
\{\widetilde{h}_n,\widetilde{k}_n\}
=
(Q|_{\mathfrak{M}})^{-1}
\bigl(Q\{\widetilde{x}_n,\widetilde{y}_n\}\bigr)
\to \{0,0\}
\quad \text{as} \quad n\to\infty.
\]
Consequently,
\[
\{u_n,\widetilde{v}_n\}
=
\{\widetilde{x}_n,\widetilde{y}_n\}
-
\{\widetilde{h}_n,\widetilde{k}_n\}
\to
\{f,\widetilde{g}\}
\quad \text{as} \quad n\to\infty.
\]

For every \(n\in\mathbb N\), the membership
\(\{u_n,\widetilde{v}_n\}\in R|_{\mathfrak H}\) implies that
\(\{u_n,\widetilde{P}_{\mathfrak H}\widetilde{v}_n\}\in\CH(R)\).
Since \(u_n\in\mathfrak H\) for every \(n\in\mathbb N\) and
\(f\in\mathfrak H\), continuity of the orthogonal projection
onto \(\mathfrak H^2\) in
\(\bigl(\widetilde{\mathfrak H}^2, \dahip\bigr)\) yields
\[
\{u_n,\widetilde{P}_{\mathfrak H}\widetilde{v}_n\}
\to
\{f,\widetilde{P}_{\mathfrak H}\widetilde{g}\}
\quad \text{as} \quad n\to\infty.
\]
Thus
\(\{f,\widetilde{P}_{\mathfrak H}\widetilde{g}\}
\in\overline{\CH(R)}\).
Since \(\{f,\widetilde{g}\}\in\widebar{R}\) with
\(f\in\mathfrak H\) was arbitrary, this proves that
\(\CH(\widebar{R})\subseteq\overline{\CH(R)}\).
Together with the reverse inclusion, this proves the lemma.
\end{proof}

\begin{proof}[Proof of Theorem~{\rm\ref{equivalences}}]
We establish the implications in the following diagram:
\[
\begin{gathered}
\text{(\ref{equivalences-i1})}
\overset{\textbf{1}}{\Rightarrow} 
\begin{gathered}[t]
\text{(\ref{equivalences-i1a})}
\\[-1mm]
\mathord{\Updownarrow}\rlap{\hspace{0.1em}\textbf{\scriptsize 2}}
\\[-1mm]
\text{(\ref{equivalences-i1b})}
\end{gathered}
\overset{\textbf{1}}{\Rightarrow} \text{(\ref{equivalences-i2})}
\overset{\textbf{3}}{\Rightarrow} \text{(\ref{equivalences-i5})}
\overset{\textbf{9}}{\Rightarrow} \text{(\ref{equivalences-i1})},
\\[1mm]
\text{(\ref{equivalences-i1})}
\overset{\textbf{4}}{\Rightarrow} \text{(\ref{equivalences-i3})}
\overset{\textbf{5}}{\Rightarrow} \text{(\ref{equivalences-i4})}
\overset{\textbf{9}}{\Rightarrow} \text{(\ref{equivalences-i1})},
\\[1mm]
\text{(\ref{equivalences-cc})}
\overset{\mkern-4mu\textbf{10}}{\Rightarrow} \text{(\ref{equivalences-i1})}
\overset{\textbf{7}}{\Rightarrow} \text{(\ref{equivalences-i20})}
\overset{\textbf{8}}{\Rightarrow} \text{(\ref{equivalences-Friedrichs})}
\overset{\mkern-4mu\textbf{10}}{\Rightarrow} \text{(\ref{equivalences-i1})} 
\overset{\textbf{6}}{\Rightarrow}\text{(\ref{equivalences-cc})}.
\end{gathered}
\]
The boldface numbers indicate the parts of the proof
in which the corresponding implications are established.
Along the indicated arrows, every statement follows
from (\ref{equivalences-i1}) and in turn implies
(\ref{equivalences-i1}).

\medskip
\noindent
\textbf{1.} (\ref{equivalences-i1})\(\Rightarrow\)(\ref{equivalences-i1a}) and
(\ref{equivalences-i1a})\(\Rightarrow\)(\ref{equivalences-i2}).
Assume that \(\widehat{\mathfrak H}\) is finite dimensional,
and let \(R\) be a closed linear relation in
\(\widetilde{\mathfrak H}\).
As a subspace of a finite-dimensional space, 
\(\widetilde{P}_{\widehat{\mathfrak H}}(\operatorname{dom} R)\)
is closed. This proves (\ref{equivalences-i1a}).  
Next, assume (\ref{equivalences-i1a}), and let \(R\) be
an arbitrary closed linear relation in
\(\widetilde{\mathfrak H}\).
Then \(\widetilde{P}_{\widehat{\mathfrak H}}(\dom R)\)
is closed, so Theorem~\ref{th-comadjn} gives
\(\CH(R^{\widetilde{*}})=\CH(R)^*\).
Since \(R\) was arbitrary, this proves
(\ref{equivalences-i2}).

\medskip
\noindent
\textbf{2.} (\ref{equivalences-i1a})\(\Leftrightarrow\)(\ref{equivalences-i1b}).
Assume (\ref{equivalences-i1a}), and let \(R\) be a closed
linear relation in \(\widetilde{\mathfrak H}\).
Since \(R^{\widetilde{*}}\) is closed, (\ref{equivalences-i1a}) applied to
\(R^{\widetilde{*}}\) implies that
\(\widetilde{P}_{\widehat{\mathfrak H}}(\operatorname{dom} R^{\widetilde{*}})\)
is closed. Applying Theorem~\ref{th-domRandPR*} to \(R^{\widetilde{*}}\),
we conclude that
\(\widetilde{P}_{\mathfrak H}R^{\widetilde{*}\widetilde{*}}\) is closed.
Since \(R^{\widetilde{*}\widetilde{*}}=R\), statement (\ref{equivalences-i1b}) follows.

Conversely, assume (\ref{equivalences-i1b}), and let \(R\) be
a closed linear relation in \(\widetilde{\mathfrak H}\).
Applying (\ref{equivalences-i1b}) to the closed linear relation \(R^{\widetilde{*}}\)
shows that \(\widetilde{P}_{\mathfrak H}R^{\widetilde{*}}\) is closed.
Theorem~\ref{th-domRandPR*} then implies that
\(\widetilde{P}_{\widehat{\mathfrak H}}(\operatorname{dom} R)\)
is closed. This proves (\ref{equivalences-i1a}).

\medskip
\noindent
\textbf{3.} 
(\ref{equivalences-i2})\(\Rightarrow\)(\ref{equivalences-i5}).
Let \(R\) be a self-adjoint relation in
\(\widetilde{\mathfrak H}\).
Then \(R\) is closed, so (\ref{equivalences-i2}) gives
\(\CH(R)=\CH(R^{\widetilde{*}})=\CH(R)^*\).
Thus \(\CH(R)\) is self-adjoint in \(\mathfrak H\), proving (\ref{equivalences-i5}). 

\medskip
\noindent
\textbf{4.} 
(\ref{equivalences-i1})\(\Rightarrow\)(\ref{equivalences-i3}).
Assume that \(\widehat{\mathfrak H}\) is finite dimensional,
and let \(R\) be a maximal dissipative linear relation in
\(\widetilde{\mathfrak H}\).
Since
\(\widetilde{P}_{\widehat{\mathfrak H}}(\operatorname{dom} R^{\widetilde{*}})\)
is finite dimensional and therefore closed,
Theorem~\ref{th-comadj4}(\ref{th-comadj4-i2}) implies that
\(\CH(R)\) is maximal dissipative in \(\mathfrak H\), proving (\ref{equivalences-i3}). 

\medskip
\noindent
\textbf{5.} 
(\ref{equivalences-i3})\(\Rightarrow\)(\ref{equivalences-i4}).
Assume (\ref{equivalences-i3}), and let \(R\) be a maximal
symmetric linear relation in \(\widetilde{\mathfrak H}\).
By Theorem~\ref{th-maxnnR}(\ref{th-maxnnR-i2}), either \(R\)
or \(-R\) is maximal dissipative.
Statement (\ref{equivalences-i3}) shows that either \(\CH(R)\) or
\(\CH(-R)=-\CH(R)\) is maximal dissipative.
Moreover, \(\CH(R)\) is symmetric by
Proposition~\ref{pr-inher}(\ref{pr-inher-i2}).
Another application of
Theorem~\ref{th-maxnnR}(\ref{th-maxnnR-i2}) therefore shows
that \(\CH(R)\) is maximal symmetric.

\medskip
\noindent
\textbf{6.} 
(\ref{equivalences-i1})\(\Rightarrow\)(\ref{equivalences-cc}).
This implication is proved in Lemma~\ref{le-comcl}.

\medskip
\noindent
\textbf{7.} 
(\ref{equivalences-i1})\(\Rightarrow\)(\ref{equivalences-i20}).
Assume that \(\widehat{\mathfrak H}\) is finite dimensional 
and let \(R\) be an arbitrary linear relation in
\(\widetilde{\mathfrak H}\).
We have already proved that both
(\ref{equivalences-cc}) and (\ref{equivalences-i2}) hold
under this assumption.
Applying (\ref{equivalences-i2}) to the closed linear relation
\(\widebar{R}\), and using
\((\widebar{R})^{\widetilde{*}}=R^{\widetilde{*}}\), gives 
\(\CH(R^{\widetilde{*}})=\CH(\widebar{R})^*\). 
By (\ref{equivalences-cc}), 
\[
\CH(R^{\widetilde{*}})
=\CH(\widebar{R})^*
=\bigl(\overline{\CH(R)}\bigr)^*
=\CH(R)^*.
\]
Since \(R\) was arbitrary, (\ref{equivalences-i20}) is proved.

\medskip
\noindent 
\textbf{8.} 
(\ref{equivalences-i20})\(\Rightarrow\)(\ref{equivalences-Friedrichs}).
Let \(S\) be a nonnegative linear relation in \(\widetilde{\mathfrak H}\).
By Proposition~\ref{pr-inhnn}, the compression \(\CH(S)\) is nonnegative in \(\mathfrak H\). Statement (\ref{equivalences-i20}) gives
\(\CH(S^{\widetilde{*}})=\CH(S)^*\).
Theorem~\ref{th-FriCom} therefore implies that
\(\CH(S_{\mathrm F})=(\CH(S))_{\mathrm F}\), proving
(\ref{equivalences-Friedrichs}).

\medskip
\noindent
\textbf{9.} 
(\ref{equivalences-i4}) and (\ref{equivalences-i5}) each imply (\ref{equivalences-i1}).
We prove these two implications by contraposition.
Assume that \(\widehat{\mathfrak H}\) is infinite dimensional.
Since \(\mathfrak H\neq\{0\}\) and its inner product is
nondegenerate, there exists \(h_0\in\mathfrak H\) such that
\([h_0,h_0]_{\mathfrak H}\neq0\).
Set \(\mathfrak H_0=\operatorname{span}\{h_0\}\), and let
\(\mathfrak H_1\) be its orthogonal complement in \(\bigl(\mathfrak{H},\kip_{\mathfrak{H}} \bigr)\).
Then \(\mathfrak H_0\) and \(\mathfrak H_1\) are Krein
subspaces and
\(\mathfrak H=\mathfrak H_0[\dotplus]\mathfrak H_1\).
The Krein subspace \(\mathfrak K = \mathfrak H_0[\dotplus]\widehat{\mathfrak H}\) 
of \(\widetilde{\mathfrak H}\) is infinite dimensional.
By Lemma~\ref{le-nidom}, there exists a densely defined
self-adjoint operator \(B_0\) in \(\mathfrak K\) such that
\(h_0\notin\operatorname{dom} B_0\).
With respect to the orthogonal decomposition
\(\widetilde{\mathfrak H}=\mathfrak K[\dotplus]\mathfrak H_1\),
define \(B=B_0\oplus0_{\mathfrak H_1}\). 
Then \(B\) is a densely defined self-adjoint operator in
\(\widetilde{\mathfrak H}\).

We first compute its compression to \(\mathfrak H\). We claim \((\operatorname{dom} B)\cap\mathfrak H=\mathfrak H_1\). The inclusion \(\supseteq\) is clear. To prove the converse let \(h\in(\operatorname{dom} B)\cap\mathfrak H\), and write
\(h=\alpha h_0+h_1\), where \(\alpha\in\mathbb C\) and
\(h_1\in\mathfrak H_1\).
Since \(\mathfrak H_1\subseteq\operatorname{dom} B\),
we have \(\alpha h_0=h-h_1\in\operatorname{dom} B\). The fact that \(h_0\notin\operatorname{dom} B\) forces
\(\alpha=0\), proving the claim. Also, \(Bh_1=0\) for all \(h_1\in\mathfrak H_1\). The claim implies 
\[
\CH(B)=\mathfrak H_1\times\{0\}.
\]

The linear relation \(\mathfrak H_1\times\mathfrak H_0\) is symmetric:
the two subspaces are orthogonal in \(\mathfrak H\).
Since \(\{0,h_0\}\) belongs to this relation but does not
belong to \(\CH(B)\), it is a proper symmetric extension
of \(\CH(B)\). Thus \(\CH(B)\) is not maximal symmetric and therefore is
not self-adjoint. On the other hand, \(B\) is self-adjoint and hence maximal
symmetric. It follows that (\ref{equivalences-i5}) and  (\ref{equivalences-i4}) respectively fail.

\medskip
\noindent
\textbf{10.} 
(\ref{equivalences-cc}) and
(\ref{equivalences-Friedrichs}) each imply (\ref{equivalences-i1}). 
Again, we proceed by contraposition.
Assume that \(\widehat{\mathfrak H}\) is infinite dimensional,
and let \(B\) be the operator constructed in part~\textbf{9}.
Set
\[
S=(\operatorname{dom} B)\times\{0\}.
\]
Since \(\operatorname{dom} B\) is dense in
\(\widetilde{\mathfrak H}\), we have
\(\widebar{S}=\widetilde{\mathfrak H}\times\{0\}\).
Moreover,
\(\CH(S)=\mathfrak H_1\times\{0\}\), which is closed.
Therefore,
\[
\CH(\widebar{S})=\mathfrak H\times\{0\},
\qquad
\overline{\CH(S)}=\mathfrak H_1\times\{0\}.
\]
These two linear relations are different, because
\(h_0\in\mathfrak H\setminus\mathfrak H_1\).
Thus (\ref{equivalences-cc}) fails.

Finally, \(S\) is nonnegative, because
\([0,f]_{\widetilde{\mathfrak H}}=0\) for all
\(f\in\operatorname{dom} B\).
Its adjoint is
\(S^{\widetilde{*}}=\widetilde{\mathfrak H}\times\{0\}\).
Every vector in \(\widetilde{\mathfrak H}\) is the limit of a sequence of vectors in \(\operatorname{dom} B\), and the inner product condition in Definition~\ref{de-Fe} is identically zero for pairs in \(S\). That definition therefore gives \(S_{\mathrm F}=\widetilde{\mathfrak H}\times\{0\}\). 

On the other hand, \(\CH(S)=\mathfrak H_1\times\{0\}\), and its adjoint in \(\mathfrak H\) is \(\mathfrak H\times\mathfrak H_0\). Since \(\mathfrak H_1\) is closed, the approximation condition in Definition~\ref{de-Fe} forces the first component of a
pair in \((\CH(S))_{\mathrm F}\) to belong to \(\mathfrak H_1\).

Conversely, for every \(u\in\mathfrak H_1\) and every \(v\in\mathfrak H_0\), the pair \(\{u,v\}\) belongs to \(\CH(S)^*\), and the constant sequence of pairs \(\{u,0\}\in\CH(S)\) satisfies the approximation and the inner product conditions. Hence
\[
(\CH(S))_{\mathrm F}
=\mathfrak H_1\times\mathfrak H_0,
\qquad
\CH(S_{\mathrm F})=\mathfrak H\times\{0\}.
\]
Since \(h_0\neq0\), the pair \(\{0,h_0\}\) belongs to
the first relation but not to the second.
Thus (\ref{equivalences-Friedrichs}) also fails.

The proof is complete.
\end{proof}

\begin{remark}
Notice that (\ref{equivalences-i2}) is equivalent to: For every linear relation \(R\) in \(\widetilde{\mathfrak H}\) the linear relation \(\CH(\widebar{R})\) is closed and \(\CH(\widebar{R}) = \CH(R^{\widetilde{*}})^*\). 
\end{remark}

\appendix

\section{Functionals} \label{Sec-funs} \label{appendix}

In the appendix we characterize the adjoint identity \(\CH(R^{\widetilde{*}})=\CH(R)^*\) and the closedness of  \(\widetilde{P}_{\mathfrak{H}}R^{\widetilde{*}}\) in terms of boundedness conditions on linear functionals on \(R\) in \(\widetilde{\mathfrak H}\).
We use Riesz functionals \(\psi_{\widetilde{u}, \widetilde{v}}\) in the Krein space \(\bigl( \widetilde{\mathfrak H}^2, \lfkip \bigr)\) defined for \(\{\widetilde{u}, \widetilde{v}\} \in \widetilde{\mathfrak H}^2\) by
\[
\psi_{\widetilde{u}, \widetilde{v}}\bigl(\{\widetilde{f},\widetilde{g}\}\bigr)
=
\bigl\llf
\{\widetilde{f},\widetilde{g}\},
\{\widetilde{u}, \widetilde{v}\}
\bigr\rlf,
\qquad
\{\widetilde{f},\widetilde{g}\}\in R.
\]

\begin{proposition} \label{pr-sfs}
Let \(R\) be a linear relation in the Krein space \(\widetilde{\mathfrak H}=\mathfrak H [\dotplus] \widehat{\mathfrak H}\).  
For all \(\{\widetilde{u}, v\} \in \widetilde{\mathfrak H}\times\mathfrak H\), the following statements hold. 

\begin{enumerate}
\renewcommand*\theenumi{\alph{enumi}}
\renewcommand*\labelenumi{\rm(\theenumi)}

\item \label{pr-sfs-i1}
If \(R\) is closed, then \(R|_{\mathfrak H}\subseteq\nul \psi_{\widetilde u,v}\) if and only if \(\{\widetilde{u}, v\} \in \overline{\widetilde P_{\mathfrak H}R^{\widetilde{*}}}\). 

\item \label{pr-sfs-i2}
\(\{\widetilde u, v\}\in \widetilde P_{\mathfrak H}R^{\widetilde{*}}\) if and only if there 
exists a constant \(c_{\widetilde u, v}>0\) such that
\begin{equation}\label{eq-ADW-estimate2}
\bigl|
\psi_{\widetilde u, v}\bigl(\{\widetilde f,\widetilde g\}\bigr)
\bigr|
\leq
c_{\widetilde u, v} \,
\bigl\|
\widetilde P_{\widehat{\mathfrak H}}\widetilde f
\bigr\|_{\widehat{\mathfrak H}} 
\quad \text{for all} \quad 
\{\widetilde f,\widetilde g\}\in R.
\end{equation}
\end{enumerate}
\end{proposition}

\begin{proof}
Let \(\{\widetilde{u}, v\} \in \widetilde{\mathfrak H}\times\mathfrak H\) be arbitrary. 

(\ref{pr-sfs-i1})  
Assume that \(R\) is closed. For all \(\{f,\widetilde g\}\in\mathfrak H\times\widetilde{\mathfrak H}\)
and all \(\{\widetilde x,\widetilde y\}\in R^{\widetilde{*}}\) we have
\[
\bigl\llf
\{f,\widetilde g\},
\{\widetilde x,\widetilde P_{\mathfrak H}\widetilde y\}
\bigr\rlf
=
\bigl\llf
\{f,\widetilde g\},
\{\widetilde x,\widetilde y\}
\bigr\rlf. 
\]
Therefore  
\[
\bigl( \widetilde P_{\mathfrak H}R^{\widetilde{*}} \bigr)^{\llf\perp\rlf}\cap
(\mathfrak H\times\widetilde{\mathfrak H}) = (R^{\widetilde{*}})^{\llf\perp\rlf}\cap 
(\mathfrak H\times\widetilde{\mathfrak H}).  
\]
Since \(R\) is closed we have \((R^{\widetilde{*}})^{\llf\perp\rlf} = R\), and consequently  
\[
\bigl( \widetilde P_{\mathfrak H} R^{\widetilde{*}} \bigr)^{\llf\perp\rlf}\cap
\bigl(\mathfrak{H}\times\widetilde{\mathfrak H}\bigr) = R|_{\mathfrak H}. 
\]
Since \(\bigl(\mathfrak{H}\times\widetilde{\mathfrak H}\bigr)^{\llf\perp\rlf} = \{0\}\times\widehat{\mathfrak{H}}\), taking orthogonal complements in \(\bigl(\widetilde{\mathfrak{H}}^2, \lfkip\bigr)\) gives 
\[
\bigl(R|_{\mathfrak{H}}\bigr)^{\llf\perp\rlf}
= \overline{\overline{\widetilde P_{\mathfrak H}R^{\widetilde{*}}}
+ \bigl(\{0\}\times\widehat{\mathfrak{H}}\bigr)}.
\]
The subspaces
\(\overline{\widetilde P_{\mathfrak H}R^{\widetilde{*}}}\)
and \(\{0\}\times\widehat{\mathfrak{H}}\) are closed
and orthogonal in the Hilbert space
\(\bigl(\widetilde{\mathfrak{H}}^2,\dahip\bigr)\),
because the first is contained in
\(\widetilde{\mathfrak{H}}\times\mathfrak{H}\).
Their sum is therefore closed, and the preceding equality becomes
\[
\bigl(R|_{\mathfrak{H}}\bigr)^{\llf\perp\rlf}
=
\overline{\widetilde P_{\mathfrak H}R^{\widetilde{*}}}
\ \lda\dotplus\rda\
\bigl(\{0\}\times\widehat{\mathfrak{H}}\bigr).
\]
Intersecting the preceding equality with \(\widetilde{\mathfrak{H}}\times\mathfrak{H}\) yields
\[
\bigl(R|_{\mathfrak H}\bigr)^{\llf\perp\rlf} \cap \bigl(\widetilde{\mathfrak H}\times\mathfrak H\bigr) = \overline{\widetilde P_{\mathfrak H}R^{\widetilde{*}}}. 
\]
Recall that \(\{\widetilde{u}, v\} \in \widetilde{\mathfrak H}\times\mathfrak H\), and we have \(R|_{\mathfrak H}\subseteq\nul \psi_{\widetilde u,v}\) if and only if \(\{\widetilde{u}, v\}\in \bigl(R|_{\mathfrak H}\bigr)^{\llf\perp\rlf}\). Consequently, the last displayed equality proves (\ref{pr-sfs-i1}).

(\ref{pr-sfs-i2})   
Assume first that \(\{\widetilde u, v\}\in \widetilde P_{\mathfrak H}R^{\widetilde{*}}\). Then there exists
\(\widehat w\in\widehat{\mathfrak H}\) such that \(\{\widetilde u, v+\widehat w\}\in R^{\widetilde{*}}\). Therefore, for all \(\{\widetilde f,\widetilde g\}\in R\), we have \(\bigl\llf
\{\widetilde f,\widetilde g\}, \{\widetilde u,v+\widehat w\} \bigr\rlf = 0\); that is 
\[
\psi_{\widetilde u, v}\bigl(\{\widetilde f,\widetilde g\}\bigr)
=
-
\bigl\llf
\{\widetilde f,\widetilde g\},
\{0,\widehat w\}
\bigr\rlf = - \iu \bigl[\widetilde f, \widehat w \bigr] = - \iu \bigl[\widetilde P_{\widehat{\mathfrak H}} \widetilde f, \widehat w \bigr]. 
\]
By the Cauchy--Bunyakovsky--Schwarz inequality in the Hilbert space
\(\bigl(\widehat{\mathfrak H}, \ahip_{\widehat{\mathfrak H}} \bigr)\) we have  
\[
\Bigl|
\psi_{\widetilde u, v}\bigl(\{\widetilde f,\widetilde g\}\bigr)
\Bigr| = 
\Bigl|
\Bigl\langle 
\widetilde P_{\widehat{\mathfrak H}}\widetilde{f},\widehat{J} \widehat{w}
\Bigr\rangle_{\widehat{\mathfrak H}}
\Bigr|
\leq
\bigl\|
 \widehat{w} 
\bigr\|_{\widehat{\mathfrak H}}
\bigl\|
\widetilde P_{\widehat{\mathfrak H}}\widetilde f
\bigr\|_{\widehat{\mathfrak H}} 
\quad \text{for all} \quad 
\{\widetilde f,\widetilde g\}\in R.
\]
Thus \eqref{eq-ADW-estimate2} holds.

Conversely, assume \eqref{eq-ADW-estimate2}. Then \(R|_{\mathfrak H}\subseteq\nul \psi_{\widetilde u,v}\). Hence the formula
\[
\varphi\bigl(
\widetilde P_{\widehat{\mathfrak H}}\widetilde f
\bigr)
=
-\psi_{\widetilde u, v}\bigl(\{\widetilde f,\widetilde g\}\bigr),
\quad
\{\widetilde f,\widetilde g\}\in R,
\]
defines a continuous linear functional on
\(\widetilde{P}_{\widehat{\mathfrak H}}(\dom R)\). This functional can be extended by continuity to the closure \(\overline{\widetilde{P}_{\widehat{\mathfrak H}}(\dom R)}\) and by zero to the orthogonal complement of this closure in \(\bigl(\widehat{\mathfrak H},\ahip_{\widehat{\mathfrak H}}\bigr)\). By the Riesz representation theorem  in the Hilbert space \(\widehat{\mathfrak H}\), there exists
\(\widehat w\in\widehat{\mathfrak H}\) such that for all \(\{\widetilde f,\widetilde g\}\in R\) we have 
\[
0 = \psi_{\widetilde u, v}\bigl(\{\widetilde f,\widetilde g\}\bigr)
+
\bigl\llf
\{\widetilde f,\widetilde g\},
\{0,\widehat w\}
\bigr\rlf
= 
\bigl\llf
\{\widetilde f,\widetilde g\},
\{\widetilde u, v\}
\bigr\rlf + 
\bigl\llf
\{\widetilde f,\widetilde g\},
\{0,\widehat w\}
\bigr\rlf . 
\]
Hence \(\{\widetilde u, v+\widehat w\}\in R^{\widetilde{*}}\), and therefore \(\{\widetilde u, v\}\in \widetilde P_{\mathfrak H}R^{\widetilde{*}}\). 
This proves (\ref{pr-sfs-i2}).
\end{proof}

\begin{theorem} \label{th-iffCAACn}
Let \(R\) be a linear relation in the Krein space \(\widetilde{\mathfrak H}=\mathfrak H [\dotplus] \widehat{\mathfrak H}\). The following statements are equivalent.
\begin{enumerate}
\renewcommand*\theenumi{\roman{enumi}}
\renewcommand*\labelenumi{\rm(\theenumi)}

\item \label{pr-iffCAACn-i1}
\(\CH(R^{\widetilde{*}}) = \CH(R)^*\).

\item \label{pr-iffCAACn-i2}
For every \(\{u,v\} \in \CH(R)^*\) there exists a constant \(c_{u,v} \gt 0\) such that
\begin{equation*} % \label{eq-bddRuvn}
\Bigl| \bllf \{\widetilde{f},\widetilde{g}\}, \{u,v\}\brlf \Bigr|
\leq c_{u,v} \bigl\|\widetilde{P}_{\widehat{\mathfrak H}} \widetilde{f}  \bigr\|_{\widehat{\mathfrak H}} \quad \text{for all} \quad 
\{\widetilde{f},\widetilde{g}\} \in R.
\end{equation*}
\end{enumerate}
\end{theorem}
\begin{proof}
By Proposition~\ref{pr-inher}(\ref{pr-inher-i1}) the inclusion \(\CH(R^{\widetilde{*}}) \subseteq  \CH(R)^*\) holds and the inclusion \(\CH(R^{\widetilde{*}}) \subseteq  \widetilde{P}_{\mathfrak H} R^{\widetilde{*}}\) is a consequence of the definition of the compression. 

Assume \(\CH(R^{\widetilde{*}}) = \CH(R)^*\). Then \(\CH(R)^* \subseteq  \widetilde{P}_{\mathfrak H} R^{\widetilde{*}}\), and (\ref{pr-iffCAACn-i2}) follows from Proposition~\ref{pr-sfs}(\ref{pr-sfs-i2}). 

Assume (\ref{pr-iffCAACn-i2}). Then, by Proposition~\ref{pr-sfs}(\ref{pr-sfs-i2}) we have \(\CH(R)^* \subseteq  \widetilde{P}_{\mathfrak H} R^{\widetilde{*}}\). Consequently, 
\[
\CH(R)^* = 
\CH(R)^* \cap \bigl( \mathfrak{H}\times \widetilde{\mathfrak H} \bigr) 
\subseteq \bigl(\widetilde{P}_{\mathfrak H} R^{\widetilde{*}}\bigr) \cap \bigl(\mathfrak{H}\times \widetilde{\mathfrak H}\bigr) =  \CH(R^{\widetilde{*}}),  
\]
which, together with \(\CH(R^{\widetilde{*}}) \subseteq  \CH(R)^*\), proves (\ref{pr-iffCAACn-i1}). 
\end{proof}

\begin{theorem}\label{th-cpfn}
Let \(R\) be a closed linear relation in the Krein spece \(\widetilde{\mathfrak H}=\mathfrak H [\dotplus] \widehat{\mathfrak H}\). The following statements are equivalent.
\begin{enumerate}
\renewcommand*\theenumi{\roman{enumi}}
\renewcommand*\labelenumi{\rm(\theenumi)}
\item\label{pr-cpfn-i0}
\(\widetilde{P}_{\mathfrak H}R^{\widetilde{*}}\) is closed in \(\widetilde{\mathfrak H}\times\mathfrak H\). 

\item\label{pr-cpfn-i2}
For every continuous linear functional \(\phi: R \to \mathbb C\) such that \(R|_{\mathfrak H}  \subseteq \nul\phi\), there exists a constant \(c_\phi>0\) such that
\[
\Bigl|\phi\bigl(\{\widetilde f,\widetilde g\}\bigr)\Bigr|
\le c_\phi \mkern2mu
\bigl\|
\widetilde P_{\widehat{\mathfrak H}}\widetilde f
\bigr\|_{\widehat{\mathfrak H}} 
\quad \text{for all} \quad
\{\widetilde f,\widetilde g\}\in R.
\]
\end{enumerate}
\end{theorem}

\begin{proof} 
The implication (\ref{pr-cpfn-i2})\(\Rightarrow\)(\ref{pr-cpfn-i0}) follows from Proposition~\ref{pr-sfs}. Indeed, if
\(\{\widetilde u,v\}\in\overline{\widetilde{P}_{\mathfrak H}R^{\widetilde{*}}}\), then
Proposition~\ref{pr-sfs}{\rm(\ref{pr-sfs-i1})} yields \(R|_{\mathfrak H}\subseteq\nul\psi_{\widetilde u,v}\). By (\ref{pr-cpfn-i2}), the functional \(\psi_{\widetilde u,v}\) satisfies the
estimate in Proposition~\ref{pr-sfs}{\rm(\ref{pr-sfs-i2})}. Hence \(\{\widetilde u,v\}\in\widetilde P_{\mathfrak H}R^{\widetilde{*}}\). Thus \(\widetilde P_{\mathfrak H}R^{\widetilde{*}}\) is closed.

To prove the converse (\ref{pr-cpfn-i0})\(\Rightarrow\)(\ref{pr-cpfn-i2}), assume that \(\widetilde P_{\mathfrak H}R^{\widetilde{*}}\) is closed in \(\widetilde{\mathfrak H}\times\mathfrak H\). Let \(\phi:R\to\mathbb C\) be a continuous linear functional such that \(R|_{\mathfrak H}\subseteq\nul\phi\). 
Since \(R\) is closed, \(\phi\) is represented on \(R\) by a Riesz functional
on \(\widetilde{\mathfrak H}^2\). Thus there exists
\(\{\widetilde u,\widetilde v\}\in\widetilde{\mathfrak H}^2\) such that
\[
\phi\bigl(\{\widetilde f,\widetilde g\}\bigr)
=
\bigl\llf
\{\widetilde f,\widetilde g\},
\{\widetilde u,\widetilde v\}
\bigr\rlf,
\qquad
\{\widetilde f,\widetilde g\}\in R.
\]
Set
\[
v=\widetilde P_{\mathfrak H}\widetilde v,
\quad
\widehat v=\widetilde P_{\widehat{\mathfrak H}}\widetilde v.
\]
Then \(\widetilde v=v+\widehat v\), and
\[
\phi
=
\psi_{\widetilde u,v}
+
\psi_{0,\widehat v},
\]
where
\[
\psi_{0,\widehat v}\bigl(\{\widetilde f,\widetilde g\}\bigr)
=
\bigl\llf
\{\widetilde f,\widetilde g\},
\{0,\widehat v\}
\bigr\rlf = \iu \bigl[\widetilde{f},\widehat{v}\bigr]_{\widetilde{\mathfrak H}}, 
\quad
\{\widetilde f,\widetilde g\}\in R.
\]

We first prove that \(R|_{\mathfrak H}\subseteq\nul\psi_{\widetilde u,v}\). 
Let \(\{f,\widetilde g\}\in R|_{\mathfrak H}\). Then \(f\in\mathfrak H\), and
therefore
\[
\psi_{0,\widehat v}\bigl(\{f,\widetilde g\}\bigr)
=
\bigl\llf
\{f,\widetilde g\},
\{0,\widehat v\}
\bigr\rlf
=
\iu [f,\widehat v]_{\widetilde{\mathfrak H}}
=0.
\]
Since \(\phi\) vanishes on \(R|_{\mathfrak H}\), it follows that \(\psi_{\widetilde u,v}\bigl(\{f,\widetilde g\}\bigr)=0\). 
Thus \(R|_{\mathfrak H}\subseteq\nul\psi_{\widetilde u,v}\).

By Proposition~\ref{pr-sfs}{\rm(\ref{pr-sfs-i1})} we have \(\{\widetilde u,v\}\in\overline{\widetilde P_{\mathfrak H}R^{\widetilde{*}}}\).  Since \(\widetilde P_{\mathfrak H}R^{\widetilde{*}}\) is closed by assumption, we get \(\{\widetilde u,v\}\in\widetilde P_{\mathfrak H}R^{\widetilde{*}}\). Therefore, by Proposition~\ref{pr-sfs}{\rm(\ref{pr-sfs-i2})}, there exists
\(c_{\widetilde u,v}>0\) such that
\begin{equation} \label{eq-psif}
\bigl|
\psi_{\widetilde u,v}\bigl(\{\widetilde f,\widetilde g\}\bigr)
\bigr|
\leq
c_{\widetilde u,v} 
\bigl\|
\widetilde P_{\widehat{\mathfrak H}}\widetilde f
\bigr\|_{\widehat{\mathfrak H}} 
\quad \text{for all} \quad 
\{\widetilde f,\widetilde g\}\in R.
\end{equation}

On the other hand,
\[
\psi_{0,\widehat v}\bigl(\{\widetilde f,\widetilde g\}\bigr)
=
\bigl\llf
\{\widetilde f,\widetilde g\},
\{0,\widehat v\}
\bigr\rlf
=
\iu[\widetilde f,\widehat v]_{\widetilde{\mathfrak H}}
=
\iu \bigl[
\widetilde P_{\widehat{\mathfrak H}}\widetilde f,\widehat v
\bigr]_{\widehat{\mathfrak H}} \quad \text{for all} \quad 
\{\widetilde f,\widetilde g\}\in R. 
\]
By the Cauchy--Bunyakovsky--Schwarz inequality in the Hilbert space
\(\bigl(\widehat{\mathfrak H}, \ahip_{\widehat{\mathfrak H}} \bigr)\) we have  
\begin{equation} \label{eq-psihat}
\Bigl|
\psi_{0,\widehat v}\bigl(\{\widetilde f,\widetilde g\}\bigr)
\Bigr| = 
\Bigl|
\bigl\langle 
\widetilde P_{\widehat{\mathfrak H}}\widetilde{f},\widehat{J} \widehat{v}
\bigr\rangle_{\widehat{\mathfrak H}}
\Bigr|
\leq
\bigl\|
 \widehat{v} 
\bigr\|_{\widehat{\mathfrak H}}
\bigl\|
\widetilde P_{\widehat{\mathfrak H}}\widetilde f
\bigr\|_{\widehat{\mathfrak H}} 
\quad \text{for all} \quad 
\{\widetilde f,\widetilde g\}\in R.
\end{equation}
Combining the estimates \eqref{eq-psif} and \eqref{eq-psihat}, we obtain
\[
\bigl|
\phi\bigl(\{\widetilde f,\widetilde g\}\bigr)
\bigr|
\leq
\bigl(c_{\widetilde u,v} + \bigl\|
 \widehat{v} 
\bigr\|_{\widehat{\mathfrak H}} \bigr)
\bigl\|
\widetilde P_{\widehat{\mathfrak H}}\widetilde f
\bigr\|_{\widehat{\mathfrak H}}\quad \text{for all} \quad 
\{\widetilde f,\widetilde g\}\in R.
\]
Thus (\ref{pr-cpfn-i2}) holds.
\end{proof}

\section*{Declarations}

\paragraph{\textbf{Funding}}
No funds, grants, or other support was received.

\paragraph{\textbf{Competing interests}}
The authors have no relevant financial or nonfinancial interests to disclose.

\paragraph{\textbf{Data Availability}}
\noindent No datasets were generated or analyzed during this study.

\paragraph{\textbf{Acknowledgment of AI assistance}}
\noindent During the preparation of this manuscript, the authors used ChatGPT, an AI language model developed by OpenAI, as an auxiliary tool for editorial suggestions, \LaTeX formatting, and discussion of possible proof strategies.

\end{document}